\documentclass[12pt]{amsart}

\usepackage{amsthm}

\usepackage{amssymb}
\usepackage{verbatim}
\usepackage[curve,matrix,arrow]{xy}
\usepackage{amsmath,amsfonts}
\usepackage{graphicx}
\usepackage{epsf}

\usepackage{graphics}

\newcommand{\mathsym}[1]{{}}

\newcommand{\thmref}[1]{Theorem~\ref{#1}}
\newcommand{\propref}[1]{Proposition~\ref{#1}}
\newcommand{\lemref}[1]{Lemma~\ref{#1}}
\newcommand{\eqnref}[1]{Equation~(\ref{#1})}
\newcommand{\remref}[1]{Remark~\ref{#1}}
\newcommand{\corref}[1]{Corollary~\ref{#1}}

\newcommand{\figref}[1]{Figure~\ref{#1}}

\newcommand{\conjref}[1]{Conjecture~\ref{#1}}

  {\end{list}}

\def\vc{\vec{v}}

\def\li{L_{i}}
\def\ri{R_{i}}
\def\ddx{{\frac{d}{dx}}}
\def\oga{{\overline{\ga}}}

\def\NN{{\mathbb N}}

\newtheorem{theorem}{Theorem}[section]

\newtheorem{corollary}[theorem]{Corollary}
\newtheorem{conjecture}[theorem]{Conjecture}
\newtheorem{lemma}[theorem]{Lemma}
\newtheorem{proposition}[theorem]{Proposition}

\theoremstyle{example}

\newtheorem{remark}[theorem]{Remark}

\theoremstyle{definition}
\newtheorem*{definition}{Definition}
\theoremstyle{notation}

\newcommand{\dd}[1]{\delta_{#1}}
\newcommand{\DD}[1]{\Delta_{#1}}
\newcommand{\nn}[1]{\mu_{#1}}

\newcommand{\hj}[3]{\hat{j}_{#1}(#2,#3)}

\newcommand{\ga}{\Gamma}

\newcommand{\tg}{\tau(\Gamma)}
\newcommand{\ta}[1]{\tau(#1)}

\newcommand{\ee}[1]{E(#1)}
\newcommand{\vv}[1]{V(#1)}
\newcommand{\va}{\upsilon}
\newcommand{\vb}{\text{v} \hspace{0.5 mm}}

\newcommand{\gc}{g(C)}

\newcommand{\tc}{\theta}

\newcommand{\pp}{p_{i}}

\newcommand{\qq}{q_{i}}
\newcommand{\opp}{\overline{p}_i}
\newcommand{\orj}{\overline{R}_j}
\newcommand{\orcj}{\overline{R}_{{c_j,\opp}}}
\newcommand{\ov}{\overline{\upsilon}_i}

\newcommand{\mucan}{{\mu_\text{can}}}

\def\can{{\mathop{\rm can}}}

\def\CC{{\mathbb C}}

\def\cC{{\mathcal C}}

\def\BDV{{\mathop{\rm BDV}}}

\def\<{\langle }
\def\>{\rangle }
\newcommand{\secref}[1]{\S\ref{#1}}

\def\ed{\epsilon_{D}}

\def\inf{\text{inf}}

\def\KB{\overline{K}}

\def\gc{\bar{g}}
\def\ed{\epsilon(\ga)}
\def\aga{a(\ga)}
\def\cgc{c(\gc)}
\def\bq{\textbf{q}}
\def\vg{\varphi (\ga)}
\def\elg{\ell (\ga)}
\newcommand{\tcg}{\theta (\ga)}
\def\lag{\lambda (\ga)}
\def\min{\text{min}}

\def\diag{\text{diag}}
\newcommand{\am}{\mathrm{A}}

\newcommand{\dm}{\mathrm{D}}

\newcommand{\jm}{\mathrm{J}}

\newcommand{\lm}{\mathrm{L}}
\newcommand{\plm}{\mathrm{L^+}}

\begin{document}

\title[Zhang's Conjecture and The Effective Bogomolov Conjecture]
{Zhang's Conjecture and the Effective Bogomolov Conjecture over function fields}

\author{Zubeyir Cinkir}
\address{Zubeyir Cinkir\\
Wolfram Research\\
100 Trade Center \\
Champaign, Illinois 61820\\
USA}
\email{zubeyirc@wolfram.com}


\thanks{I would like to thank R. Rumely and M. Baker for their continued support and guidance. I would like to thank X. Faber for letting me know about Zhang's Conjecture and for useful correspondence during the preparation of this paper.}

\begin{abstract}
We prove the Effective Bogomolov Conjecture, and so the Bogomolov Conjecture, over a function field of characteristic $0$ by proving Zhang's Conjecture
about certain invariants of metrized graphs. In the function field case, these conjectures were previously known to be true only for curves of genus at most 4 and a few other special cases. We also either verify or improve the previous results. We relate the invariants involved in Zhang's Conjecture to the tau constant of metrized graphs. Then we use and extend
our previous results on the tau constant. By proving another Conjecture of Zhang, we obtain a new proof of the slope inequality for Faltings heights on moduli space of curves.
\end{abstract}

\maketitle

\section{Introduction}\label{section introduction}

In this paper, we study various invariants associated to a given metrized graph and polarized metrized graph. We derive formulas relating the invariants studied in the papers \cite{Zh2} and \cite{Fa} in terms of the tau constant of metrized graphs. This enables us to use the tools developed (
\cite{CR}, \cite{BRh}, \cite{BF}, \cite{C1}, \cite{C2}, \cite{C3}, and \cite{C4}) to study the tau constant. We extend our previous results on the tau constant (\cite{C1}, \cite{C2} and \cite{C3}), prove S. Zhang's Conjecture \cite[4.1.2]{Zh2}, and prove stronger version of X. Faber's Conjecture \cite[1.3]{Fa}. The consequences of these conjectures include the following applications in number theory and algebraic geometry:

$(i)$ We prove the effective Bogomolov's Conjecture over function fields of characteristic $0$. If a conjecture (see the articles \cite[1.4.1]{Zh2} and \cite{GS}) due to Grothendieck and Gillet-Soul\'{e} is true, our results extend to the function field of positive characteristic case, and have implications to number field case. The Bogomolov Conjecture over function fields were previously known only in some special cases, which will be discussed briefly in the next section.

$(ii)$ We give a new proof of a slope inequality for Faltings heights on moduli space of curves by proving another Conjecture of Zhang \cite[Conjecture 1.4.5]{Zh2}. This slope inequality was first proved by A. Moriwaki in the article \cite[Theorem D at page 3]{AM4} in the characteristic $0$ case and in the article \cite[Theorem 4.1]{AM2} for arbitrary characteristic.
Our method depends only on calculations involving invariants of metrized graphs, and makes it possible to obtain stronger versions of the slope inequality in certain special cases.

Throughout the paper, we use the interpretation of metrized graphs as resistive electric circuits and related electrical properties such as circuit reductions. Whenever it is needed, we consider metrized graphs only with their combinatorial graph structure. Our previous results on the tau constant (\cite{C1}, \cite{C2}, and \cite{C3}) use the properties of a continuous Laplacian operator on a metrized graph, which was defined by M. Baker and R. Rumely \cite{BRh} and studied from the perspective of harmonic analysis.
We will give a brief description of metrized graphs, and recall results from electric circuit theory and combinatorial graph theory with short explanations. Interested readers should consult the references cited in the related sections.

\section{The Bogomolov Conjecture, the slope inequality and main results}\label{sec Bogomolov Conjecture and Slope Inequality}

We first recall some definitions. Let $X$ be a smooth projective surface over a field $k$, and let $Y$ be a
smooth projective curve over $k$. A fibration $f:X \longrightarrow Y$ over $Y$ is called \textit{``isotrivial"}, if all smooth fibers are isomorphic to a fixed curve.

Let $k$ be a field. Let $X$ be a smooth projective surface over $k$, and let $Y$ be a smooth projective curve over $k$. Let $f:X \rightarrow Y$ be a semi-stable fibration such that the generic fiber of $f$ is smooth and of genus $\gc \geq 2$. Let $K$ be the function field of $Y$, with algebraic closure $\KB$, and let $C$ be the generic fiber of $f$.
The N\'{e}ron-Tate height pairing on the Jacobian variety $\text{Jac}(C)(\KB)=Pic^0(C)(\KB)$ of the curve $C/K$ induces a seminorm $|| \cdot ||_{NT}$. For $\frac{\omega_{C}}{2 \gc -2} \in Pic^1(C)(\KB)$, we have a canonical inclusion $j:C(\KB) \longrightarrow \text{Jac}(C)(\KB)$ defined by $j(x)=(2 g-2) x-\omega_{C}$.

If we define
$B_{C}(P;r)= \{ x \in C(\KB): \, ||j(x)-P||_{NT} \leq r \} $, where
$P \in Pic^0(C)(\KB)$ and $r \geq 0$, and if we set
$$r_C(P)=\begin{cases}
-\infty & \text{ if }  \#(B_C(P;0))=\infty, \\
\sup \{ r \geq 0 \, | \, \, \#(B_C(P;r)) <\infty \} & \text{otherwise.}
\end{cases} $$
then Bogomolov's conjecture can be stated as follows:
\begin{conjecture}\cite{AM2}(\textit{Bogomolov Conjecture})\label{Conj Bogomolov2}
If $f$ is non-isotrivial, then $r_C(P)>0$ for all $P$.
\end{conjecture}
\begin{conjecture} \cite{KY1}(\textit{Effective Bogomolov Conjecture})\label{Conj Bogomolov3}
If $f$ is non-isotrivial, then there exists an ``effectively calculable" positive number $r_0$ such that
$$ \inf_{P \in Pic^0(C)(\KB)} r_C(P) \geq r_0.$$
\end{conjecture}
We will now describe how metrized graphs can be related to above
conjectures by following the articles \cite{AM2} and \cite{Zh1}.

\begin{figure}
\centering
\includegraphics[scale=0.9]{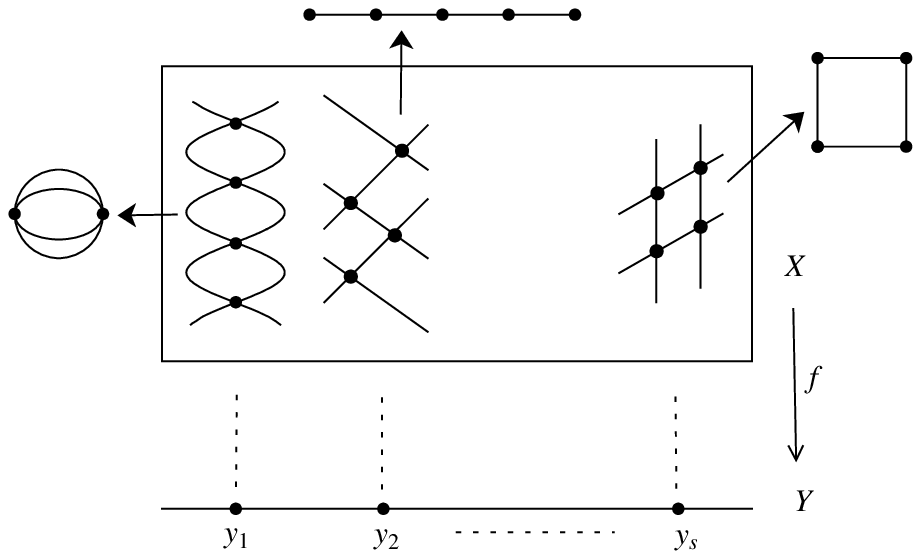} \caption{Singular fibers of $f:X\rightarrow Y$, and dual graphs of
fibers.}\label{fig Bogomolov1}
\end{figure}

For the semistable fibration $f:X \longrightarrow Y$, let $CV(f)= \{y_1,y_2, \cdots, y_s \}$ be the set of
critical values of $f$, where $s$ is the number of singular fibers. That is, $y \in CV(f)$
iff $f^{-1}(y)$ is singular. For any $y_i \in CV(f)$, let $\ga_{y_i}$ be the dual graph of the fiber $C_{y_i}:=f^{-1}(y_i)$
, for each $1 \leq i \leq s$.
The metrized graph $\ga_{y_i}$ is obtained as follows. The set of vertices $V_{y_i}$ of $\ga_{y_i}$ is indexed by irreducible
components of the fiber $f^{-1}(y_i)$ and the singularities of $f^{-1}(y_i)$ correspond to edges of length $1$. Let
$I(C_{y_i}):=\{ C_{1,y_i},C_{2,y_i}, \cdots, C_{v_i,y_i} \}$ be the set of irreducible components of the fiber $C_{y_i}$, where $v_i$
is the number of irreducible components in $C_{y_i}$. Then the irreducible curve $C_{j,y_i}$ corresponds to the
vertex $p_j \in V_{y_i}$ for each $1 \leq j \leq v_i$ (see \figref{fig dual3} and \figref{fig Bogomolov1}).
Let $\delta_{y_i}$ be the number of singularities in $C_{y_i}$.
By our construction, $\delta_{y_i}=\ell(\ga_{y_i})$, the length of $\ga_{y_i}$, for each $1 \leq i \leq s$.
Let $\delta :=\sum_{i=1}^s\delta_{y_i}$,  the total number of singularities in the fibration. For any $p_j \in V_{y_i}$,
let $\bq (p_j):=g(C_{j,y_i})$, where $g(C_{j,y_i})$ is the arithmetic genus of $C_{j,y_i}$ (see \secref{sec pm graph} for the role of $\bq$). Let $g(Y)$ be the genus of $Y$.
%

We have $K_X$, $K_Y$, $\omega_X$, $\omega_Y$, and $\omega_{X/Y}$, which are the canonical divisor of $X$, the canonical divisor of $Y$,
the dualizing sheaf of $X$, the dualizing sheaf of $Y$, and  the relative dualizing sheaf, respectively.
For the admissible dualizing sheaf $\omega_a$ \cite{Zh1} (see also the article \cite[pg. 3]{AM1}), we have the following inequalities due to Zhang \cite{Zh1}:
\begin{equation}\label{eqn admissible inequality}
\begin{split}
\omega_{X/Y}^2 \geq \omega_a^2 \geq 0.
\end{split}
\end{equation}
Zhang \cite{Zh1} showed that  $\omega_a^2 > 0$ is equivalent to the Bogomolov conjecture, i.e., \conjref{Conj Bogomolov2}.
\begin{theorem}\cite[Theorem 2.1]{AM3} \cite[Theorem 5.6]{Zh1}\label{thm inf and admis}
If $\omega_a^2 >0$, then
$$
\inf_{P \in Pic^0(C)(\KB)} r_C(P) \geq \sqrt{(g-1) \omega_a^2}.
$$
\end{theorem}
Let $\ed$ and $\vg$ be the invariants of a dual graph $\ga$ defined in \secref{sec pm graph} below.
We have (\cite[Equation 1.2.1]{Zh2}, \cite[Equation 2.2]{AM3}, and \cite[Section 4.5]{C1})
\begin{equation}\label{eqn admissible sheaves1}
\begin{split}
\omega_a^2=\omega_{X/Y}^2-\sum_{i=1}^{s}\epsilon(\Gamma_{y_i}).
\end{split}
\end{equation}

Zhang defined the canonical Gross-Schoen cycle $\Delta_{\xi}$ associated to $X$, and showed in the article \cite[Corollary 1.3.2]{Zh2} that
\begin{equation}\label{eqn admissible sheaves2}
\begin{split}
\omega_a^2=\frac{2 \gc-2}{2 \gc+1} \langle  \Delta_{\xi}, \Delta_{\xi} \rangle +\frac{2 \gc-2}{2 \gc+1} \sum_{i=1}^s\varphi (\Gamma_{y_i}).
\end{split}
\end{equation}
\begin{remark}\label{rem varphi positivity and Bogomolov Conj}
Since the height of $\Delta_{\xi}$, namely $\langle \Delta_{\xi}, \Delta_{\xi} \rangle$, is non-negative whenever the characteristic of $k$ is $0$ (as in \thmref{thm effective bogomolov conj holds} and \thmref{thm effective bogomolov conj holds 2nd embedding}), proving the positivity of $\vg$ for any polarized metrized graph (pm-graph in short) $\ga$ will be enough to prove Bogomolov Conjecture.
\end{remark}
\begin{remark}\label{rem Parsin}
Whenever $f$ is smooth, we clearly have $\omega_a^2=\omega_{X/Y}^2$ by \eqnref{eqn admissible sheaves1}, and that $\omega_{X/Y}^2\geq 12$ as Par$\check{s}$in \cite{P} showed.
\end{remark}
The Effective Bogomolov Conjecture (\conjref{Conj Bogomolov3}) was known to be true for curves of genus less than $5$ (\cite{AM1}, \cite{AM2}, \cite{AM3}, \cite{AM4}, \cite{KY1}, \cite{KY2}, and \cite{Fa}). Also, W. Gubler \cite{G} showed that the Bogomolov Conjecture is true for $C$ if the Jacobian variety of $C$ has totally degenerate reduction over some point $y \in Y$.  One can consult the article \cite{Fa} to see the list of previously known lower bounds to $\omega_a^2$.

We will use the following notation for the singularities that are in the fibers of $f$:

Let $y \in CV(f)$, and let $p \in f^{-1}(y)$ be a node. If the partial normalization of $f^{-1}(y)$ at $p$ is connected,
we say that $p$ is of type $0$. If it is disconnected, then it has two components, in which case $p$ will be said to be of type $i$,
where $i$ is the minimum of the arithmetic genera of the components. We denote the total number of nodes of type $i$ in the fiber $f^{-1}(y)$ by  $\delta_i (\Gamma_{y})$, and we set $\delta_i (X)=\sum_{j=1}^s \delta_i (\Gamma_{y_j})$. With our earlier notation, we have $\delta_{y_j}=\sum_{i \geq 0} \delta_i (\Gamma_{y_j})$ and $\delta=\sum_{i \geq 0} \delta_i (X)$.

Next we state Zhang's first conjecture which implies \conjref{Conj Bogomolov3}, and so \conjref{Conj Bogomolov2}.
\begin{conjecture}\cite[ Conjecture 1.4.2]{Zh2}\label{conj Zhang varphi}
For any $y \in CV(f)$, there is a positive continuous function  $c(\gc)$ of $\gc \geq 2$ such that the following inequality holds:
$$\varphi(\Gamma_{y}) \geq c(\gc) \delta_0(\Gamma_{y})+\sum_{i \geq 1} \frac{2 i (\gc-i)}{\gc}\delta_i(\Gamma_{y}).$$
\end{conjecture}
Previous results, due to Moriwaki and K. Yamaki when $\gc$ is $2$ or $3$, on \conjref{Conj Bogomolov3} depend on a slope inequality which is the following lower bound for $\deg f_*(\omega_{X/Y})$:
\begin{equation}\label{eqn slope inequality}
\begin{split}
\deg f_*(\omega_{X/Y}) \geq \frac{\gc}{8 \gc +4} \delta_0 (X) + \sum_{i \geq 1} \frac{i (\gc -i)}{2 \gc+1} \delta_i(X).
\end{split}
\end{equation}
Note that the inequality (\ref{eqn slope inequality}) was proved by Moriwaki \cite[Theorem D at page 3]{AM4} in the characteristic $0$ case. This slope inequality is actually closely related to the inequality given in the article \cite[Theorem 4.1]{AM2}, which is slightly weaker but holds in any characteristic. Its proof and connection to Bogomolov Conjecture is through the following equation which is obtained by Noether's formula
\begin{equation}\label{eqn Noether's formula}
\begin{split}
\deg f_*(\omega_{X/Y})=\frac{1}{12} (\omega_{X/Y}^2+\sum_{i=1}^s \delta_{y_i}).
\end{split}
\end{equation}
Let $\lambda (\Gamma)$ and $\aga$ be the invariants of a dual graph $\ga$ defined in \secref{sec pm graph} below. We have Zhang's second conjecture leading to a second proof of the slope inequality given in (\ref{eqn slope inequality}).
\begin{conjecture}\cite[ Conjecture 1.4.5]{Zh2}\label{conj Zhang lambda}
For any $y \in CV(f)$, the following inequality holds:
$$\lambda (\Gamma_{y}) \geq \frac{\gc}{8 \gc +4} \delta_0(\Gamma_{y})+\sum_{i \geq 1} \frac{ i (\gc-i)}{2 \gc +1}\delta_i(\Gamma_{y}).$$
\end{conjecture}
Zhang reduced and unified his first and second conjectures into the following conjecture:
\begin{conjecture}\cite[Conjecture 4.1.2]{Zh2}\label{Conj Zhang's second}
Let $\ga$ be an irreducible polarized metrized graph of genus $\gc$. Then the following two inequalities hold:
\begin{equation*}
\begin{split}
\frac{\gc-1}{\gc+1}(\ell(\ga)-4 \gc \cdot \aga) \leq \ed \leq 12 \gc \cdot \aga - (1+\cgc) \ell(\ga),
\end{split}
\end{equation*}
where $\cgc$ is a positive number for each $\gc \geq 2$.
\end{conjecture}
Faber verified this conjecture for curves of genus less than $5$.
In the rest of this section, we will state our main results.

We prove that \conjref{conj Zhang varphi} holds as follows:
\begin{theorem}\label{thm Zhang's first conjecture holds}
Let $\ga$ be a pm-graph with genus $\gc$. Then we have
\begin{equation*}
\begin{split}
\vg & \geq t(\gc) \delta_0(\ga) +\sum_{i \geq 1} \frac{2 i (\gc-i)}{\gc}\delta_i(\ga).
\end{split}
\end{equation*}
where $t(2)=\frac{1}{27}$, $t(3)=\frac{892 - 11 \sqrt{79}}{14580} \approx 0.054473927$, and $t(\gc)=\frac{(\gc-1)^2}{2 \gc (7 \gc +5)}$ for $\gc \geq 4$. In particular, $t(\gc)\geq \frac{3}{88}$ for $\gc \geq 4$.
\end{theorem}
\begin{proof}
The result follows from \thmref{thm proof of Zhangs varphi conjecture}, \thmref{thm g=3} and the article \cite[Corollary 4.4.2]{Zh2}.
\end{proof}
We believe that the lower bounds in \thmref{thm Zhang's first conjecture holds} can be improved when $\gc \geq 3$. In the light of \propref{prop varphi for complete graph}, which gives an exact formula for $\vg$ when $\ga$ is a complete graph,  the lower bounds in \thmref{thm Zhang's first conjecture holds} are not far from optimal. In \secref{sec pm graph with zero part}, we give better bounds to $\vg$ for certain classes of pm-graphs. \corref{cor varphi necklace} shows that $\vg$ can be very large for some pm-graphs.

\begin{figure}
\centering
\includegraphics[scale=0.87]{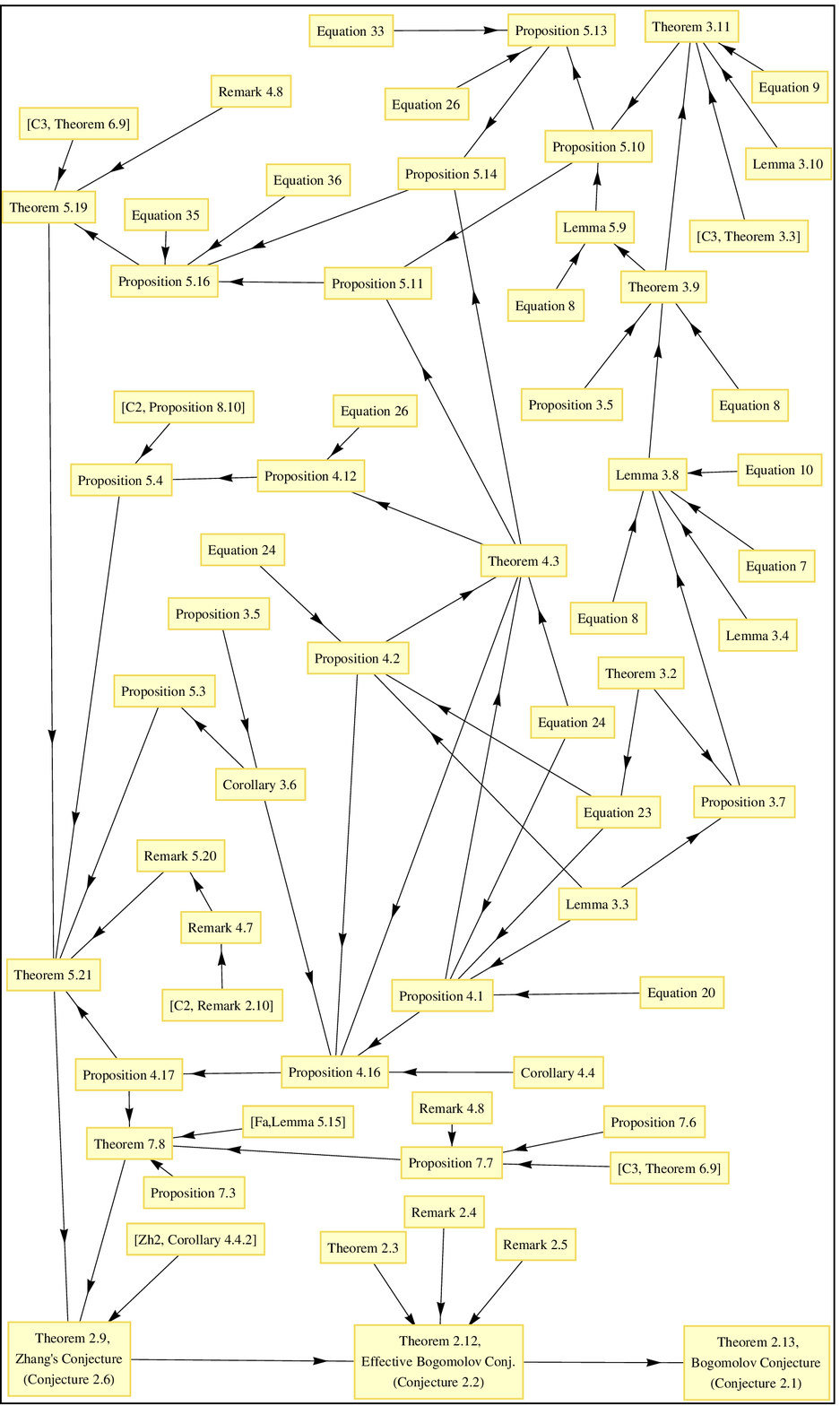} \caption{Major dependencies between the results that lead to the proofs of Conjectures (\ref{conj Zhang varphi}), (\ref{Conj Bogomolov3}), and (\ref{Conj Bogomolov2}).} \label{fig Map}
\end{figure}

We prove that \conjref{conj Zhang lambda} holds as follows (with notation as in \secref{sec metrized graphs}):
\begin{theorem}\label{thm Zhang's second conjecture holds}
Let $(\ga,\bq)$ be a pm-graph.
Then we have
\begin{equation*}
\begin{split}
\lag \geq \frac{\gc}{8 \gc +4} \delta_0 (\ga)+\sum_{i \geq 1} \frac{ i (\gc-i)}{2 \gc +1}\delta_i(\ga).
\end{split}
\end{equation*}
\end{theorem}
\begin{proof}
The result follows from \propref{prop lambda} and the article \cite[Corollary 4.4.2]{Zh2}.
\end{proof}
For any given genus $\gc$, we have examples showing that $\lag$ can be very close to the lower bound given in \thmref{thm Zhang's second conjecture holds} or it can be linear in $\gc$ (see \corref{cor varphi necklace}). Therefore, the topology of the pm-graph $\ga$ plays important role in the value of $\lag$. The techniques we have developed can be used for numeric calculations of $\lag$ in general and symbolic calculations of $\lag$ in specific cases.
\begin{theorem}\label{thm Zhang's conjecture holds}
\conjref{Conj Zhang's second} holds with $\cgc=4 t(\gc)$ for $\gc \geq 2$, where $t(\gc)$ is as in \thmref{thm Zhang's first conjecture holds}.
\end{theorem}
\begin{proof}
The result follows from \thmref{thm Zhang's first conjecture holds}, \thmref{thm Zhang's second conjecture holds} and \lemref{lem Zhang's second conj upper bound}.
\end{proof}
\begin{theorem}\label{thm effective bogomolov conj holds}
Let \text{char}$(k)=0$. If $f$ is non-isotrivial, then we have $ \inf_{P \in Pic^0(C)(\KB)} r_C(P) \geq \sqrt{r_0}$, where $r_0$ can be taken as follows:
\begin{equation*}\label{}
\begin{split}
r_0=\begin{cases}12 (\gc -1), & \text{if $f$ is smooth}
\\\frac{2 (\gc -1)^2}{2 \gc +1} \Big ( t(\gc)\delta_0(\ga)  +\sum_{i \geq 1} \frac{2 i (\gc-i)}{\gc}\delta_i(\ga) \Big ), & \text{otherwise}
\end{cases}
\end{split}
\end{equation*}
with $t(\gc)$ as in \thmref{thm Zhang's first conjecture holds}. Therefore, the Effective Bogomolov Conjecture holds.
\end{theorem}
\begin{proof}
If $f$ is smooth, then the result follows from \remref{rem Parsin} and \thmref{thm inf and admis}. Suppose that
$f$ is not smooth, then $f$ has places of bad reduction, i.e., $CV(f)$ is non-empty. Then the result follows from
\remref{rem varphi positivity and Bogomolov Conj}, \thmref{thm Zhang's first conjecture holds}, and \thmref{thm inf and admis}.
\end{proof}
Then \conjref{Conj Bogomolov2} follows from \thmref{thm effective bogomolov conj holds}:
\begin{theorem}\label{thmcor Bogomolov conj holds}
The Bogomolov Conjecture holds if \text{char}$(k)=0$.
\end{theorem}
The map of the results that lead to \thmref{thmcor Bogomolov conj holds}, \thmref{thm effective bogomolov conj holds} and \thmref{thm Zhang's first conjecture holds} can be found in \figref{fig Map}.

If we consider the slightly different embedding $j_D:C(\KB) \longrightarrow \text{Jac}(C)(\KB)$ defined by $j_D(x)= x-D$ for a given
$D \in \text{Div}^1(C(\KB))$, the Bogomolov Conjecture concerns the positivity of $\text{inf}_{D \in Div^1(C(\KB))} a'(D)$ when $f$ is non-isotrivial, where $a'(D):=\lim_{x \in C(\KB)} \, \text{inf} \, \hat{h}(j_D(x))$ and $\hat{h}$ is the
N\'{e}ron-Tate height pairing on the Jacobian variety $\text{Jac}(C)(\KB)=Pic^0(C)(\KB)$ of the curve $C/K$.
With this embedding, we can state \thmref{thm effective bogomolov conj holds} as follows:
\begin{theorem}\label{thm effective bogomolov conj holds 2nd embedding}
Let \text{char}$(k)=0$. If $f$ is non-isotrivial, then we have
\begin{equation*}\label{}
\begin{split}
\inf_{D \in \text{Div}^1(C(\KB))} a'(D) \geq \begin{cases}\frac{3}{\gc -1}, & \text{if $f$ is smooth}
 \\\frac{1}{2(2 \gc +1)} \Big ( t(\gc)\delta_0(\ga)  +\sum_{i \geq 1} \frac{2 i (\gc-i)}{\gc}\delta_i(\ga) \Big ), & \text{otherwise}
\end{cases}
\end{split}
\end{equation*}
with $t(\gc)$ as in \thmref{thm Zhang's first conjecture holds}. Therefore, the Effective Bogomolov Conjecture holds.
\end{theorem}
\thmref{thm effective bogomolov conj holds 2nd embedding} shows that a stronger version of Faber's Conjecture \cite[1.3]{Fa} holds.

Note that \thmref{thm effective bogomolov conj holds}, \thmref{thmcor Bogomolov conj holds} and \thmref{thm effective bogomolov conj holds 2nd embedding} will also extend to the positive characteristic case if the following conjecture holds:
\begin{conjecture}(Grothendieck, Gillet-Soul\'{e}) \cite{GS} \label{conj positivity of Gross-Shoen cycle}
Let $k$ be a number field or a function field with positive characteristic, then the following height inequality holds:
$$\langle  \Delta_{\xi}, \Delta_{\xi} \rangle \geq 0.$$
Moreover, it becomes an equality precisely when $\Delta_{\xi}$ is rationally equivalent to $0$.
\end{conjecture}


\section{Metrized graphs and their tau constants}\label{sec metrized graphs}

Rumely introduced metrized graphs to
study arithmetic properties of curves and developed arithmetic capacity theory.
T. Chinburg and Rumely \cite{CR} used metrized graphs when they introduced their ``capacity pairing''.
Another pairing satisfying ``desirable'' properties is the ``admissible pairing on curves'' introduced by Zhang \cite{Zh1}.
Metrized graphs were used as a non-archimedean analogue of a Riemann surface
(\cite{RumelyBook}, \cite{CR}, and \cite{Zh1}).
Following Zhang's approach, A. Moriwaki used metrized graphs
and Green's functions to prove specific cases of Bogomolov's conjecture over function fields
in a series of papers, \cite{AM1}, \cite{AM2}, and \cite{AM3}. Extending Moriwaki's approach,
Yamaki \cite{KY1} proved very special cases of effective generalized Bogomolov's conjecture over function fields.

Metrized graphs arise as dual graphs of curves. Chinburg and Rumely \cite{CR} worked with a
canonical measure $\mu_{can}$ of total mass $1$ on a metrized graph $\ga$.
Similarly, Zhang worked with a measure $\mu_{ad}$ of total mass $1$ on $\ga$.
The measure $\mu_{ad}$ defined in \secref{sec pm graph} is a generalization of $\mu_{can}$ defined in this section.

A \textit{metrized graph} $\ga$ is a finite connected graph
equipped with a distinguished parametrization of each of its edges.
In particular, $\ga$ is a one-dimensional manifold except at
finitely many ``branch points''. See also the articles \cite{RumelyBook} and  \cite{Zh1}.

A metrized graph can have multiple edges and self-loops. For any given $p \in \ga$, the number of directions emanating from $p$ will be called the \textit{valence} of $p$, and will be denoted by $\va(p)$. By definition, there can be only finitely many $p \in \ga$ with $\va(p)\not=2$.

Given a metrized graph $\ga$, we will denote its set of vertices by $\vv{\ga}$.
We require that $\vv{\ga}$ is non-empty and that $p \in \vv{\ga}$ for each $p \in \ga$ with
$\va(p)\not=2$. For a given metrized graph $\ga$, it is possible to enlarge the
vertex set $\vv{\ga}$ by considering arbitrarily many valence $2$ points as vertices.

For a given metrized graph $\ga$ with vertex set $\vv{\ga}$, the set of edges of $\ga$ is the set of closed line segments with end points in $\vv{\ga}$. We will denote the set of edges of $\ga$ by $\ee{\ga}$. However, if
$e_i$ is an edge, by $\ga-e_i$ we mean the graph obtained by deleting the {\em interior} of $e_i$.

Let $v:=\# (\vv{\ga})$ and $e:=\# (\ee{\ga})$. We define the \textit{genus} of $\ga$ to be the first Betti number $g(\ga):=e-v+1$ of the graph $\ga$. We will simply use $g$ to show $g(\ga)$ when there is no danger of confusion. Note that the genus is a topological invariant of $\ga$. In particular, it is independent of the choice of the vertex set $\vv{\ga}$.
Since $\ga$ is connected, $g(\ga)$ coincides with the cyclotomic number of $\ga$ in combinatorial graph theory.

We denote the length of an edge $e_i \in \ee{\ga}$ by $\li$. Then \textit{total length} of $\ga$, which will be denoted by $\elg$, is given by $\elg=\sum_{i=1}^e\li$.

The minimum number of vertices whose deletion
disconnects $\ga$ is called the \textit{``vertex connectivity"} of $\ga$ and
will be denoted by $\kappa(\ga)$. The minimum number of edges whose deletion
disconnects $\ga$ is called the \textit{``edge connectivity"} of $\ga$ and
will be denoted by $\Lambda(\ga)$. Let $\overline{\delta}(\ga):= \min \{ \va(p) | p \in \vv{\ga} \}$ be the minimum
of valences of the vertices. Then by basic graph theory \cite[pg. 3]{BB1},
$\kappa(\ga) \leq \Lambda(\ga) \leq \overline{\delta}(\ga)$.
We call a metrized graph $\ga$ \textit{irreducible}, as in the article \cite{Fa},
if it can not be disconnected by deleting any single point. That is,
$\ga$ has vertex connectivity at least $2$ for each possible choice of
vertex set $\vv{\ga}$. Therefore, if $\ga$ is irreducible, it has edge
connectivity at least $2$. If the edge connectivity of a metrized graph $\ga$ is at least two, we also say that $\ga$ is a \text{bridgeless} metrized graph. It is clear from the definitions that every irreducible graph is bridgeless, but there can be bridgeless graphs which have vertex connectivity $1$ and so are not irreducible. For example,
union of two copies of the circle graph along a vertex is a bridgeless metrized graph but not an irreducible metrized graph.

Baker and Rumely \cite{BRh} defined the following measure valued \textit{Laplacian} on a given metrized graph:
\begin{equation}
\DD{x}(f(x))=-f''(x)dx - \sum_{p \in \vv{\ga}}\bigg[ \sum_{\vc
\hspace{0.5 mm} \text{at} \hspace{0.5 mm} p}
d_{\vc}f(p)\bigg]\dd{p}(x),
\end{equation}
for a continuous function $f : \Gamma \rightarrow \CC$ such that $f$ is
$\cC^2$ on $\ga
\backslash \vv{\ga}$ and $f^{\prime \prime}(x) \in L^1(\Gamma)$. See the article \cite[Section 4]{BRh}, for the description of $\BDV(\Gamma)$, the largest set of continuous functions for which $\DD{x}$ is defined.

In the article \cite{CR}, a kernel $j_{z}(x,y)$ giving
a fundamental solution of the Laplacian is defined and studied as a
function of $x, y, z \in \Gamma$. For fixed $z$ and $y$ it
has the following physical interpretation: when $\Gamma$ is viewed
as a resistive electric circuit with terminals at $z$ and $y$,
with the resistance in each edge given by its length, then
$j_{z}(x,y)$ is the voltage difference between $x$ and $z$,
when unit current enters at $y$ and exits at $z$ (with reference
voltage $0$ at $z$).
\begin{lemma}\cite[Lemma 2.10]{CR}\label{lemjzero} The function $j_\zeta(x,y)$ is symmetric in $x$ and in $y$,
is jointly continuous as a function of all three variables, and is
nonnegative, with $j_{\zeta}(\zeta,y) = j_{\zeta}(x,\zeta) = 0$ for
all $x, y, \zeta \in \Gamma$.
\end{lemma}
A self-contained proof of this fact is given in the article \cite{Zh1}.

The effective resistance between two points $x,y$ of a metrized graph $\ga$ is given by $r(x,y)=j_y(x,x).$
We call $j_{z}(x,y)$ and $r(x,y)$ be the \textit{voltage function} and the \textit{resistance function} on $\ga$, respectively.
The functions $j_{z}(x,y)$ and $r(x,y)$ are also studied in the articles \cite{BRh}, \cite{BF}, \cite{C1}, \cite{C2}, \cite{C3}, and \cite{C4}.

We will denote by $R_{i}(\ga)$, or by $R_i$ if there is no danger of confusion, the resistance between the end points of an edge $e_i$ of a graph $\ga$ when the interior of the edge $e_i$ is deleted from $\ga$.

\begin{figure}
\centering
\includegraphics[scale=0.7]{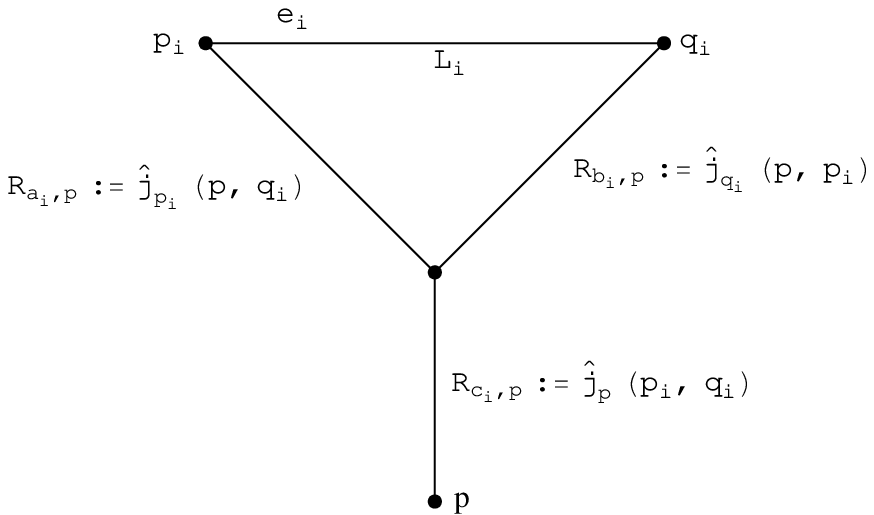} \caption{Circuit reduction with reference to an edge and a point.} \label{fig 2termp}
\end{figure}

Let $\ga$ be a metrized graph with $p \in \vv{\ga}$, and let $e_i \in \ee{\ga}$ having end points $\pp$ and $\qq$.
If $\ga -e_i$ is connected, then $\ga$
can be transformed to the graph in Figure \ref{fig 2termp}
by circuit reductions. More details on this fact can be found in the articles \cite{CR} and \cite[Section 2]{C2}.
Note that in \figref{fig 2termp}, we have $R_{a_i,p} = \hj{\pp}{p}{\qq}$,
$R_{b_i,p} = \hj{\qq}{p}{\pp}$, $R_{c_i,p} = \hj{p}{\pp}{\qq}$, where $\hj{x}{y}{z}$
is the voltage function in $\ga-e_i$. We have $R_{a_i,p}+R_{b_i,p}=\ri$ for each $p \in \ga$.

If $\ga-e_i$ is not connected, we set $R_{b_i,p}=\ri=\infty$ and $R_{a_i,p}=0$ if $p$ belongs to the component of $\ga-e_i$
containing $\pp$, and we set $R_{a_i,p}=\ri=\infty$ and $R_{b_i,p}=0$ if $p$ belongs to the component of $\ga-e_i$
containing $\qq$. We will use these notation in the rest of the paper.

For any real-valued, signed Borel measure $\mu$ on $\Gamma$ with
$\mu(\Gamma)=1$ and $|\mu|(\Gamma) < \infty$, define the function
$j_{\mu}(x,y) \ = \ \int_{\Gamma} j_{\zeta}(x,y) \, d\mu({\zeta}).$
Clearly $j_{\mu}(x,y)$ is symmetric, and is jointly continuous in
$x$ and $y$. Chinburg and Rumely \cite{CR} discovered that there is a unique real-valued, signed Borel measure $\mu=\mu_{can}$
such that $j_{\mu}(x,x)$ is constant on $\ga$. The measure $\mu_\can$ is called the
\textit{canonical measure}. See the articles \cite{BRh} and \cite{C2} for several interpretations of $\mu_{can}$.
Baker and Rumely \cite{BRh} called the constant $\frac{1}{2}j_{\mu}(x,x)$ the \textit{tau constant} of $\ga$
and denoted it by $\tg$.

Let $\mu$ be a real-valued signed Borel measure of total mass $1$ on
$\Gamma$. In the article \cite{BRh}, the \textit{Arakelov-Green's function} $g_\mu(x,y)$ associated to
$\mu$ is defined to be $g_\mu(x,y) = \int_\Gamma j_\zeta(x,y) \, d\mu(\zeta) - \int_\Gamma
j_\zeta(x,y) \, d\mu(\zeta) d\mu(x) d\mu(y)$, where the latter integral is a constant that depends on $\ga$ and $\mu$.

As shown in the article \cite{BRh}, $g_{\mu}(x,y)$ is continuous,
symmetric (i.e., $g_{\mu}(x,y)=g_{\mu}(y,x)$, for each $x$ and $y$),
and for each $y$, $\int_{\Gamma} g_{\mu}(x,y) \, d\mu(x) \ = \ 0 \ $.
More precisely, as shown in the article \cite{BRh}, one can characterize
$g_\mu(x,y)$ as the unique function on $\Gamma \times \Gamma$ such
that
\begin{itemize}
\item[$(1)$]
$g_\mu(x,y)$ is jointly continuous in $x, y$ and belongs to
$\BDV_\mu(\Gamma)$ as a
 function of $x$, for each fixed $y$, where $\BDV_\mu(\Gamma) := \{ f \in \BDV(\Gamma) \; : \; \int_\Gamma f \,
d\mu = 0 \}$.
\item[$(2)$]
For fixed $y$, $g_\mu$ satisfies the identity
$
\Delta_x g_\mu(x,y) = \delta_y(x) - \mu(x).
$
\item[$(3)$]
$\iint_{\Gamma \times \Gamma} g_\mu(x,y) d\mu(x) d\mu(y) = 0$.
\end{itemize}

The diagonal values $g_{\mu_{can}}(x,x)$ are constant on $\ga$, and are equal to the tau constant $\tg$. In terms of spectral theory, when $\ga$ has total length $1$, the tau constant is the trace of the inverse operator of $\Delta$. Note that the notation $\tg$ is used in the article \cite[Equation 4.1.2]{Zh2} to denote another invariant of $\ga$.

The following theorem gives an explicit description of $\mu_{can}$:
\begin{theorem}\cite[Theorem 2.11]{CR} \label{thmCanonicalMeasureFormula}
For a given metrized graph $\ga$, let $L_i$ be the length of edge $e_i \in \ee{\ga}$, and let $R_i$ be the effective resistance between the endpoints of $e_i$ in the graph $\Gamma-e_i$. Then we have
\begin{equation*}
\mu_\can(x) \ = \ \sum_{p \in \vv{\ga}} (1 - \frac{1}{2}\vb(p))
\, \delta_p(x) + \sum_{e_i \in \ee{\ga}} \frac{dx}{L_i+R_i}.
\end{equation*}
\end{theorem}
%
The following two lemmas express $\tg$ in terms of the resistance function and the canonical measure.
\begin{lemma}\cite{REU}\label{lemtau1}
For any metrized graph $\ga$ and its resistance function $r(x,y)$, and for each $x \in \ga$,
$\tau(\Gamma) = \frac{1}{2}\int_{\ga} r(x,y) d\mucan(y).$
\end{lemma}
\begin{lemma}\cite[Lemma 14.4]{BRh}\label{lemtauformula}
For any fixed $y \in \Gamma$, we have
$
\tau(\Gamma) = \frac{1}{4} \int_\Gamma
\left( \ddx r(x,y) \right)^2 dx .
$
\end{lemma}
\lemref{lemtauformula} implies that $\tg \geq 0$ for any metrized graph $\ga$.

For the resistance function $r(x,y)$ in $\ga$, we use circuit reductions (parallel and series reductions, see the article 
\cite[Section 2]{C2} and the related references given therein) to obtain the following equalities:
%
%
%
%
%
\begin{equation}\label{eqn2term0}
\begin{split}
r(p_i,p)=\frac{(\li+R_{b_i,p})R_{a_i,p}}{\li+\ri}+R_{c_i,p}, \quad
\text{and} \quad r(q_i,p)=\frac{(\li+R_{a_i,p})R_{b_i,p}}{\li+\ri}+R_{c_i,p}.
\end{split}
\end{equation}
Therefore,
\begin{equation}\label{eqn2term0a}
\begin{split}
r(\pp,p)-r(\qq,p) &= \frac{\li(R_{a_i,p}-R_{b_i,p})}{\li+\ri},
\\ r(\pp,p)+r(\qq,p) & = \frac{\li \ri}{\li+\ri}+2 \frac{R_{a_i,p} R_{b_i,p}}{\li+\ri}+2R_{c_i,p}.
\end{split}
\end{equation}
The following proposition is obtained by evaluating the integral formula for the tau constant, given in \lemref{lemtauformula}, on each edges of $\ga$.
\begin{proposition}\cite{REU}\label{proptau}
Let $\Gamma$ be a metrized graph, and let $L_i$ be the length of the
edge $e_{i}$, for $i \in \{1,2, \dots, e\}$.
Using the notation above,
if we fix a vertex $p$ we have
\[
\ta{\ga} = \frac{1}{12} \sum_{e_i \in \ee{\ga}} \frac{\li^3+3\li(R_{a_{i},p}-R_{b_{i},p})^2}{(\li+\ri)^2}.
\]
Here, if $\ga-e_i$ is not connected, i.e. $\ri$ is infinite, the
summand corresponding to $e_i$ should be replaced by $3\li$, its limit as $\ri \longrightarrow \infty$.
\end{proposition}
The proof of \propref{proptau} can be found in the article \cite[Proposition 2.9]{C2}.
\begin{corollary}\label{cor tau for circle}
Let $\ga$ be a circle graph. Then $\tg =\frac{1}{12} \elg$.
\end{corollary}
\begin{proof}
The result follows from \propref{proptau} (equivalently, see the article \cite[Corollary 2.17]{C2}).
\end{proof}
Chinburg and Rumely \cite[page 26]{CR} showed that
\begin{equation}\label{eqn genusa}
\sum_{e_i \in \ee{\ga}}\frac{\li}{\li +\ri}=g, \quad \text{equivalently } \sum_{e_i \in \ee{\ga}}\frac{\ri}{\li +\ri}=v-1.
\end{equation}
The following proposition gives another formula for the tau constant which depends on the expression of $\mu_\can(x)$ given in \thmref{thmCanonicalMeasureFormula}.
\begin{proposition}\label{prop tau canonical}
Let $\ga$ be a metrized graph and let $r(x,y)$ be the resistance function in $\ga$.
Then for any $p \in \vv{\ga}$,
$$ \tg = -\frac{1}{4} \sum_{q \in  \vv{\ga}}(\va(q)-2)r(p,q)
+ \frac{1}{2} \sum_{e_i \in \, \ee{\ga}} \frac{1}{\li+\ri} \int_{0}^{\li} r(p,x) dx.$$
\end{proposition}
\begin{proof} We have
\begin{equation*}
\begin{split}
& \tg = \frac{1}{2}\int_{\ga} r(p,x) d\mucan(x), \quad \text{by \lemref{lemtau1}, so
by \thmref{thmCanonicalMeasureFormula}},
\\  & = \frac{1}{2}\sum_{q \in \vv{\ga} } (1 - \frac{1}{2}\va(p)) \int_{\ga} r(p,x) \delta_q(x)
 + \frac{1}{2} \sum_{e_i \in \, \ee{\ga}} \frac{1}{\li+\ri} \int_{0}^{\li} r(p,x) dx .
\end{split}
\end{equation*}
Here, each edge $e_i \in \ee{\ga}$
is parametrized by a segment $[0,\li]$, under its arclength parametrization. Thus the result follows.
\end{proof}
The purpose of the following lemma is to clarify the relation between the formulas for $\tg$ given in \propref{proptau} and \propref{prop tau canonical}. This will help us to derive \thmref{thm tau formula new}.
\begin{lemma}\label{lem crtterm}
Let $\ga$ be a bridgeless metrized graph, and let $\pp$ and $\qq$ be the end points of $e_i \in \ee{\ga}$. For any $p
\in \vv{\ga}$, we have
\begin{equation*}
\begin{split}
\sum_{e_i \in \,
\ee{\ga}}\frac{\li(R_{a_{i},p}-R_{b_{i},p})^2}{(\li+\ri)^2} &
=\sum_{e_i \in \, \ee{\ga}}\frac{\li}{\li + \ri}
\big(r(\pp,p)+r(\qq,p)\big) - \sum_{q \in
\vv{\ga}}(\va(q)-2)r(p,q)\\
& = 2\sum_{q \in \vv{\ga}}r(p,q) -\sum_{e_i \in \,
\ee{\ga}}\frac{\ri}{\li + \ri} \big(r(\pp,p)+r(\qq,p)\big).
\end{split}
\end{equation*}
\end{lemma}
\begin{proof}
Let $p \in \vv{\ga}$. By \lemref{lemtauformula}, $ 4 \tau(\Gamma) = \int_\Gamma
\left( \ddx r(p,x) \right)^2 dx.$ Thus, integration by parts gives
\begin{equation}\label{eqn crtterm1}
\begin{split}
4 \tg &=\sum_{e_i \in \, \ee{\ga}} \big(r(p,x) \cdot \ddx
r(p,x)\big)|^{\li}_{0} - \sum_{e_i \in \, \ee{\ga}} \int_{0}^{\li}
r(p,x) \frac{d^2}{dx^2}r(p,x) dx.
\end{split}
\end{equation}
If $x \in e_i$, then by parallel and series circuit reductions applied to the graph in \figref{fig 2termp}
\begin{equation}\label{eqn crtterm resistance}
\begin{split}
r(p,x)=\frac{(\li-x+R_{b_{i},p}) (x+R_{a_{i},p})}{\li +\ri} + R_{c_{i},p}, \text{  so } \frac{d}{dx}r(p,x)=\frac{-2 x+\li+R_{b_{i},p}-R_{a_{i},p}}{\li+\ri}.
\end{split}
\end{equation}
Thus, $\frac{d^2}{dx^2}r(p,x)=\frac{-2}{\li+\ri}$ if $x \in e_i$. This equality along with
\propref{prop tau canonical}, Equations (\ref{eqn crtterm1}), (\ref{eqn2term0}), and (\ref{eqn2term0a})
give the first equality in the theorem. Then the second
equality holds by the following identity:
\begin{equation}\label{eqn crtterm2}
\begin{split}
\sum_{q \in \vv{\ga}} \va(q) r(p,q) = \sum_{e_i \in \, \ee{\ga}}
\big( r(\pp,p)+r(\qq,p)\big).
\end{split}
\end{equation}
\end{proof}
Now we are ready to state a new formula for $\tg$ which will play a crucial role in proving \conjref{Conj Zhang's second} and so the related conjectures.
\begin{theorem}\label{thm tau formula new}
Let $\ga$ be a bridgeless metrized graph, and let $r(x,y)$ be the resistance function on it. Then for any given $p \in \vv{\ga}$, we have
\begin{equation*}
\begin{split}
\tg=\frac{\elg}{12}-\frac{1}{6}\sum_{q \in \vv{\ga}} (\va(q)-2) r(p,q)+\frac{1}{3}\sum_{e_i \in \, \ee{\ga}} \frac{\li}{\li+\ri}R_{c_i,p}.
\end{split}
\end{equation*}
\end{theorem}
\begin{proof}
By \eqnref{eqn2term0a} and the first equality in \lemref{lem crtterm}, we have
\begin{equation}\label{eqn new form1}
\begin{split}
\sum_{e_i \in \,
\ee{\ga}}\frac{\li(R_{a_{i},p}-R_{b_{i},p})^2}{(\li+\ri)^2} &
=\sum_{e_i \in \, \ee{\ga}}\frac{\li^2 \ri +2 \li R_{a_{i},p} R_{b_{i},p}}{(\li + \ri)^2}
+ 2 \sum_{e_i \in \, \ee{\ga}}\frac{\li}{\li + \ri} R_{c_i,p}
\\ & \qquad - \sum_{q \in \vv{\ga}}(\va(q)-2)r(p,q).
\end{split}
\end{equation}
By using $\ri=R_{a_{i},p}+R_{b_{i},p}$, we obtain $2 R_{a_{i},p} R_{b_{i},p}=\frac{1}{2}\ri^2-\frac{1}{2}(R_{a_{i},p} -R_{b_{i},p})^2$; then substituting this in \eqnref{eqn new form1} we obtain the following equality
\begin{equation}\label{eqn new form2}
\begin{split}
3\sum_{e_i \in \,
\ee{\ga}}\frac{\li(R_{a_{i},p}-R_{b_{i},p})^2}{(\li+\ri)^2} &
=\sum_{e_i \in \, \ee{\ga}}\frac{2 \li^2 \ri +\li \ri^2}{(\li + \ri)^2}
+ 4 \sum_{e_i \in \, \ee{\ga}}\frac{\li}{\li + \ri} R_{c_i,p}
\\ & \qquad - 2 \sum_{q \in \vv{\ga}}(\va(q)-2)r(p,q).
\end{split}
\end{equation}
Adding $\sum_{e_i \in \, \ee{\ga}}\frac{\li^3}{(\li + \ri)^2}$ to both sides of \eqnref{eqn new form2} and using the fact that
$\elg =\sum_{e_i \in \, \ee{\ga}}\li$ gives
\begin{equation*}
\begin{split}
\sum_{e_i \in \,
\ee{\ga}}\frac{\li^3 + 3 \li(R_{a_{i},p}-R_{b_{i},p})^2}{(\li+\ri)^2} &
= \elg + 4 \sum_{e_i \in \, \ee{\ga}}\frac{\li}{\li + \ri} R_{c_i,p}
- 2 \sum_{q \in \vv{\ga}}(\va(q)-2)r(p,q).
\end{split}
\end{equation*}

On the other hand the left hand side is equal to $12 \tg$ by \propref{proptau}. Then the result follows.
\end{proof}

For each $i \in \{1,2, \dots e\}$, let
$\oga_i$ be the metrized graph obtained from a metrized graph $\ga$ by contracting the i-th edge
$e_{i} \in \ee{\ga}$ to its end
points. If $e_{i} \in \ga$ has end points $\pp$ and $\qq$, then in
$\oga_i$, these points become an identical vertex which we will denote as $\opp$. Note that $\ell(\ga)=\ell(\oga_i)+\li$.
We will denote the valence of $p \in \oga_i$ by $\ov(p)$. If $q \in \vv{\ga}-\{ \pp,\qq \}$, then $q \in \vv{\oga_i}$ and $\ov(q)=\va(q)$. If $e_i$ is not a self loop, then $\ov(\opp)=\va(\pp)+\va(\qq)-2$.
\begin{lemma}\label{lem contraction middle}
Let $\ga$ be a metrized graph with an edge $e_i \in \ee{\ga}$ such that $\ga-e_i$ is connected, and let $\oga_i$ be defined as before. Let $r(x,y)$ and $r_i(x,y)$ be the resistance functions in $\ga$ and $\oga_i$, respectively. Then we have
\begin{equation*}
\begin{split}
\sum_{q \in \vv{\oga_i}} (\ov(q)-2) r_i(\opp,q)=\sum_{q \in \vv{\ga}}(\va(q)-2) \Big( \frac{R_{a_{i},q} R_{b_{i},q}}{\ri} + R_{c_i,q} \Big).
\end{split}
\end{equation*}
\end{lemma}
\begin{proof}
When $e_i$ is contracted, $\ga$ shown in \figref{fig 2termp} becomes $\oga_i$ shown in
\figref{fig gammacontri}. Therefore, by applying parallel circuit reduction and using the fact that $R_{a_{i},q}+ R_{b_{i},q}=\ri$ for each $q \in \vv{\ga}$ we obtain the following equality:
\begin{equation}\label{eqn resistance for contraction graph}
\begin{split}
r_i(\opp,q) = \frac{R_{a_{i},q} R_{b_{i},q}}{\ri} + R_{c_i,q}.
\end{split}
\end{equation}

\eqnref{eqn resistance for contraction graph} and the fact that $\ov(q)=\va(q)$ for each $q \in \vv{\ga}-\{ \pp,\qq \}$ yields the following equality for each $q \in \vv{\ga}-\{ \pp,\qq \}$:

\begin{equation}\label{eqn resistance for contraction graph2}
\begin{split}
(\ov(q)-2)r_i(\opp,q) = (\va(q)-2) \Big(\frac{R_{a_{i},q} R_{b_{i},q}}{\ri} + R_{c_i,q}\Big).
\end{split}
\end{equation}

On the other hand, $R_{a_{i},q} \cdot R_{b_{i},q}=0=R_{c_i,q}$ for $q \in \{ \pp,\qq \}$, and
$r_i(\opp,\opp)=0$. Thus the result follows from \eqnref{eqn resistance for contraction graph2}.

\begin{figure}
\centering
\includegraphics[scale=0.9]{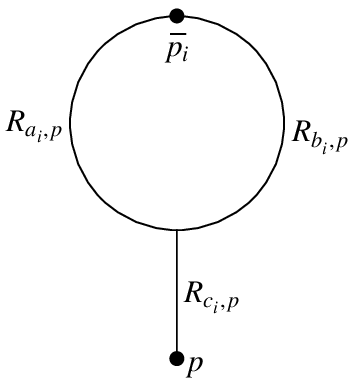}\caption{$\oga_i$ after circuit reductions.} \label{fig gammacontri}
\end{figure}
\end{proof}
Next, by using \lemref{lem contraction middle} and several other results about $\tg$ we will prove the following theorem which will be useful in proving the second inequality in \conjref{Conj Zhang's second}.
\begin{theorem}\label{thm tau new formula with contraction}
Let $\ga$ be a bridgeless metrized graph with $\# (\vv{\ga})=v \geq 3$. Let $R_{a_{i},q}$, $R_{b_{i},q}$, $\ri$ and $R_{c_i,q}$ be as defined before for each $e_i \in \ee{\ga}$, and let $\orj$ and $\orcj$ be defined similarly for each $e_j \in \ee{\oga_i}$. Then we have
\begin{equation*}
\begin{split}
\tg & = \frac{\ell (\ga)}{12} -\frac{1}{6 (v-2)} \sum_{q \in \vv{\ga}}  (\va(q)-2) \sum_{e_i \in \ee{\ga}}  \frac{R_{a_{i},q} R_{b_{i},q}+\ri R_{c_i,q}}{\li+\ri}
\\ & \qquad + \frac{1}{3 (v-2)} \sum_{e_i \in \ee{\ga}} \frac{\ri}{\li +\ri} \sum_{e_j \in \ee{\oga_i}}  \frac{L_j \orcj}{L_j+\orj}.
\end{split}
\end{equation*}
\end{theorem}
\begin{proof}
Since $\ga$ is bridgeless, $\oga_i$ is also bridgeless for each $e_i \in \ee{\ga}$. Thus, we can apply \thmref{thm tau formula new} for each $\oga_i$ with vertex $\opp$ to obtain,
\begin{equation}\label{eqn tau new for contraction}
\begin{split}
\tau(\oga_i) =\frac{\ell (\oga_i)}{12}-\frac{1}{6}\sum_{q \in \vv{\oga_i}} (\ov(q)-2) r_i(\opp,q)+\frac{1}{3}\sum_{e_j \in \ee{\oga_i}}  \frac{L_j \orcj}{L_j+\orj},
\end{split}
\end{equation}
where $r_i(\opp,q)$ is the resistance, in $\oga_i$, between the vertices $\opp$ and $q$. Then we multiply both sides of \eqnref{eqn tau new for contraction} by $\frac{\ri}{\li + \ri}$ and sum over all edges $e_i \in \ee{\ga}$ to obtain
\begin{equation}\label{eqn tau new sum for contraction}
\begin{split}
\sum_{e_i \in \ee{\ga}} \frac{\ri \tau(\oga_i)}{\li +\ri} & = \frac{1}{12} \sum_{e_i \in \ee{\ga}} \frac{\ri \ell (\oga_i)}{\li +\ri} -\frac{1}{6}\sum_{e_i \in \ee{\ga}} \frac{\ri}{\li +\ri} \sum_{q \in \vv{\oga_i}} (\ov(q)-2) r_i(\opp,q)
\\ & \qquad +\frac{1}{3}\sum_{e_i \in \ee{\ga}} \frac{\ri}{\li +\ri} \sum_{e_j \in \ee{\oga_i}}  \frac{L_j \orcj}{L_j+\orj}.
\end{split}
\end{equation}
On the other hand, it follows from
the article \cite[Theorem 3.3]{C3} that
\begin{equation}\label{eqn tau sum contraction formula}
\begin{split}
\sum_{e_i \in \ee{\ga}} \frac{\ri \tau(\oga_i)}{\li +\ri} = (v-2) \tg + \frac{1}{12} \sum_{e_i \in \ee{\ga}} \frac{\li^2}{\li +\ri}.
\end{split}
\end{equation}
Multiply both sides of the equation given in \lemref{lem contraction middle} by $\frac{\ri}{\li + \ri}$ and sum over all edges $e_i \in \ee{\ga}$ to obtain
\begin{equation}\label{eqn contraction sum middle}
\begin{split}
\sum_{e_i \in \ee{\ga}} \frac{\ri}{\li +\ri}\sum_{q \in \vv{\oga_i}} (\ov(q)-2) r_i(\opp,q)=\sum_{q \in \vv{\ga}}(\va(q)-2) \sum_{e_i \in \ee{\ga}} \frac{R_{a_{i},q} R_{b_{i},q} + \ri R_{c_i,q}}{\li +\ri}.
\end{split}
\end{equation}
The result follows by substituting \eqnref{eqn tau sum contraction formula} and \eqnref{eqn contraction sum middle} into \eqnref{eqn tau new sum for contraction}, using the fact that $\ell (\oga_i)=\elg - \li$, and using \eqnref{eqn genusa}.
\end{proof}

\section{Polarized metrized graphs}\label{sec pm graph}



\begin{figure}
\centering
\includegraphics[scale=0.85]{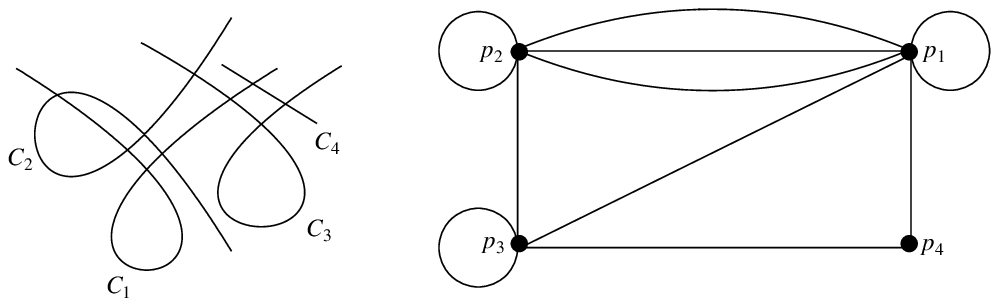} \caption{A graph and its dual
graph} \label{fig dual3}
\end{figure}

In this section, following the articles \cite{Zh2} and \cite{Fa}, we first introduce the notion of a polarized metrized graph and related concepts. Then we give formulas for several invariants $\ed$, $\vg$, $a(\ga)$, and $\lag$. These invariants will be important in proving \conjref{Conj Zhang's second}. Finally, we prove the first inequality in \conjref{Conj Zhang's second} and derive a formula for $\vg$ which will be used in \secref{sec pm graph with zero part} to prove the second inequality in \conjref{Conj Zhang's second}.

Let $\ga$ be a metrized graph and let $\bq : \ga \rightarrow  \NN$ be a function on the set of vertices of $\ga$.
The canonical divisor $K$ of $(\ga,\bq)$ is defined to be the following divisor on $\ga$:
\begin{equation}\label{eqn app2a}
\begin{split}
K  = \sum_{p \in \vv{\ga}} (\va(p)-2+2 \bq (p))p, \quad \text{and} \quad
 \dd{K}(x)  = \sum_{p \in \vv{\ga}} (\va(p)-2+\bq (p))\dd{p}(x).
\end{split}
\end{equation}
The pair $(\ga,\bq)$ will be called a \textit{polarized metrized graph} (pm-graph in short) if $\bq$ is non-negative and $K$ is an effective divisor. The \textit{genus} $\gc (\ga)$ of a pm-graph $(\ga,\bq)$ is defined to be
\begin{equation}\label{eqn genus}
\begin{split}
\gc (\ga) = 1+\frac{1}{2}\deg{K}=g(\ga)+\sum_{p \in \vv{\ga}}\bq (p).
\end{split}
\end{equation}
We will simply use $\gc$ to show $\gc(\ga)$ when there is no danger of confusion. Note that $\gc \geq 1$ for a pm-graph.
We call a pm-graph $(\ga,\bq)$ \textit{irreducible} if the underlying metrized graph $\ga$ is irreducible.

Note that the reduction graph $R(X)$ of any semistable curve $X$ of genus $\gc$ over a discrete valuation ring
is a pm-graph of genus $\gc$.

Recall that how pm-graphs are obtained from fibres of a semistable fibration $f:X \longrightarrow Y$ is explained at page $2$.

Let $\nn{ad}(x)$ be the admissible metric associated to $K$ (as defined in the article \cite[Lemma 3.7]{Zh1}). Then
\begin{equation}\label{eqn app2}
\begin{split}
\nn{ad}(x) = \frac{1}{\gc} \Big(\sum_{p \in \vv{\ga}}\bq (p) \dd{p}(x) +
\sum_{i \in \ee{\ga}} \frac{dx}{L_{i}+R_{i}} \Big).
\end{split}
\end{equation}
Then by \thmref{thmCanonicalMeasureFormula}, we have
\begin{equation}\label{eqn app2b}
\begin{split}
\nn{ad}(x) = \frac{1}{2 \gc}(2\nn{can}(x)+\dd{K}(x)).
\end{split}
\end{equation}
Moreover, $\dd{K}(\ga) =  \deg(K)  =  2 \gc-2$, and by \eqnref{eqn genusa} $\nn{can}(\ga) =  1 = \nn{ad}(\ga)$.

Set $\tcg:=\sum_{p, \, q \in \, \vv{\ga}}(\va(p)-2+2 \bq (p))(\va(q)-2+2 \bq (q))r(p,q)$, and define
\begin{equation}\label{eqn app3}
\begin{split}
 & \ed = \iint_{\Gamma \times \Gamma} r(x,y) \dd{K}(x)  \nn{ad}(x),
\\ & a(\ga) = \frac{1}{2} \iint_{\Gamma \times \Gamma} r(x,y) \nn{ad}(x) \nn{ad}(y),
\\ & \vg = 3 \gc \cdot a(\ga) -\frac{1}{4} (\ed +\ell (\ga)),
\\ & \lag = \frac{\gc-1}{6 (2 \gc+1)} \vg +\frac{1}{12}(\ed +\elg).
\end{split}
\end{equation}
We have $\tcg \geq 0$ for any pm-graph $\ga$, since the corresponding canonical divisor $K$ is effective.

Note that the second invariant in \eqnref{eqn app3} was denoted as $\tg$ in the article \cite{Zh2}. In order not to have notational conflict with the articles \cite{BRh}, \cite{C1}, \cite{C2}, \cite{C3}, and \cite{C4}, we denote it by $a(\ga)$.
For a pm-graph $(\ga,\bq)$ by $\tg$ we mean the tau constant of the underlying metrized graph $\ga$.
\begin{proposition}\label{prop ed and tau}
Let $\ga$ be a pm-graph. Then we have
\begin{equation*}
\begin{split}
\ed = \frac{(4 \gc-4) \tg}{\gc} + \frac{\tcg}{2 \gc}.
\end{split}
\end{equation*}
\end{proposition}
\begin{proof} By \eqnref{eqn app3}, we have
\begin{equation*}
\begin{split}
\ed & = \iint_{\Gamma \times \Gamma} r(x,y) \dd{K}(x)  \nn{ad}(x).  \quad
\text{Then by \eqnref{eqn app2a},}
\\ & = \sum_{p \in \,
\vv{\ga}}(\va(p)-2+2 \bq (p)) \int_{\ga} r(p,y) \nn{ad}(y), \quad \text{and by \eqnref{eqn app2b}}
\\ & = \sum_{p \in \,
\vv{\ga}}(\va(p)-2+2 \bq (p)) \int_{\ga} r(p,y) \Big( \frac{1}{2 \gc}(2\nn{can}(y)+\dd{K}(y) \Big).
\end{split}
\end{equation*}
Then the result follows by \lemref{lemtau1} and the fact that $\deg (K)=2 \gc -2$.
\end{proof}

\begin{figure}
\centerline{\epsffile{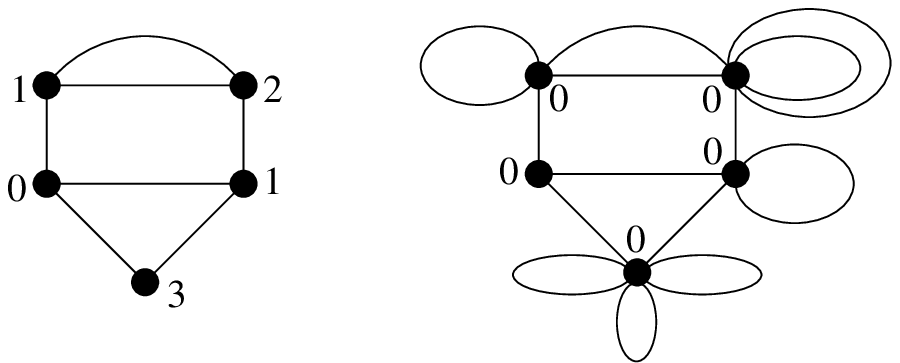}} \caption{$\ga$ and $\ga_0$. Values of $\bq$ and $\bq_0$ are shown in each graph.} \label{fig simplified}
\end{figure}

Recall that both $\tg$ and $\tcg$ are nonnegative for any pm-graph $\ga$. Therefore, \propref{prop ed and tau} implies that
$\ed \geq 0$. Similarly, $a(\ga) \geq 0$ for any pm-graph $\ga$ by the following proposition.
\begin{proposition}\label{prop adm int of res}
Let $\ga$ be a pm-graph. Then we have
\begin{equation*}
\begin{split}
a(\ga) = \frac{(2 \gc -1)\tg}{\gc^2}+\frac{\tc(\ga)}{8 \gc^2}.
\end{split}
\end{equation*}
\end{proposition}
\begin{proof} By \eqnref{eqn app3}, we have
\begin{equation*}
\begin{split}
2 a(\ga) & = \iint_{\Gamma \times \Gamma} r(x,y) \nn{ad}(x) \nn{ad}(y).  \quad
\text{Then by \eqnref{eqn app2b},}
\\ & = \int_{\ga} \Big( \int_{\ga} r(x,y) \frac{1}{2 \gc}(2\nn{can}(x)+\dd{K}(x)) \Big) \nn{ad}(y)
\\ & = \int_{\ga} \Big( \frac{2 \tg}{\gc}+\frac{1}{2 \gc}\sum_{p \in \,
\vv{\ga}}(\va(p)-2+2 \bq (p))r(p,y) \Big) \nn{ad}(y), \quad
\text{by \lemref{lemtau1}.}
\\ & = \frac{2 \tg}{\gc} +\frac{1}{2 \gc}\sum_{p \in \,
\vv{\ga}}(\va(p)-2+2 \bq (p)) \int_{\ga} r(p,y) \nn{ad}(y), \quad
\text{since $\nn{ad}(\ga)=1$}.
\end{split}
\end{equation*}
Since $\deg (K)=2 \gc -2$, the result follows by a similar calculation done to obtain the third equality from the first equality.
\end{proof}
\begin{theorem}\label{thm main1}
Let $\ga$ be a pm-graph. Then we have
\begin{equation*}
\begin{split}
\vg = \frac{(5 \gc -2) \tg}{\gc}+\frac{\tcg}{4 \gc}-\frac{\elg}{4}.
\end{split}
\end{equation*}
\end{theorem}
\begin{proof}By \eqnref{eqn app3}, we have
\begin{equation*}
\begin{split}
\vg & = 3 \gc (\ga) \cdot a(\ga) -\frac{1}{4} (\ed + \ell (\ga)). \quad
\text{Then by \propref{prop adm int of res} and \propref{prop ed and tau}}
\\ & = 3 \gc \Big(  \frac{(2 \gc -1) \tg}{\gc^2}+\frac{\tc(\ga)}{8 \gc^2} \Big) -
\frac{1}{4} \Big(\frac{(4 \gc -4) \tg}{\gc}+\frac{\tcg}{2 \gc}+\elg \Big).
\end{split}
\end{equation*}
This is equivalent to what we wanted.
\end{proof}
\begin{corollary}\label{cor lambda interms of tau}
Let $\ga$ be a pm-graph. Then we have
\begin{equation*}
\begin{split}
\lag = \frac{(3\gc -3) \tg}{4 \gc+2}+\frac{\tcg}{16 \gc+8}+\frac{(\gc +1)\elg}{16 \gc+8}.
\end{split}
\end{equation*}
\end{corollary}
\begin{proof}
By substituting the formula for $\vg$ given in \thmref{thm main1} and the formula for $\ed$ given in \propref{prop ed and tau}
into formula for $\lag$ given in \eqnref{eqn app3}, we obtain the result.
\end{proof}
\begin{proposition}\label{prop relative dual sheaf interms of tau}
Let $X$, $Y$, $\Delta_{\xi}$ and related notation be as in \secref{sec Bogomolov Conjecture and Slope Inequality}. Then,
\begin{equation*}\label{}
\begin{split}
\omega_{X/Y}^2= \frac{2 \gc-2}{2 \gc+1} \langle  \Delta_{\xi}, \Delta_{\xi} \rangle + \frac{18 (\gc -1)}{2 \gc +1} \sum_{i=1}^{s}\ta{\Gamma_{y_i}}
+\frac{3}{4 \gc +2} \sum_{i=1}^{s} \theta (\Gamma_{y_i}) -\frac{\gc -1}{4 \gc +2} \delta.
\end{split}
\end{equation*}
\end{proposition}
\begin{proof}
By \eqnref{eqn admissible sheaves1} and \eqnref{eqn admissible sheaves2},
\begin{equation}\label{eqn relative dual sheaf interms of varphi}
\begin{split}
\omega_{X/Y}^2= \frac{2 \gc-2}{2 \gc+1} \langle  \Delta_{\xi}, \Delta_{\xi} \rangle +\frac{2 \gc-2}{2 \gc+1} \sum_{i=1}^s\varphi (\Gamma_{y_i})
+\sum_{i=1}^{s}\epsilon(\Gamma_{y_i}).
\end{split}
\end{equation}
Then the result follows from \thmref{thm main1}, \propref{prop ed and tau} and the fact that $\delta=\sum_{i=1}^{s} \ell(\Gamma_{y_i})$ (see page $2$).
\end{proof}
\begin{remark}\label{rem upper bounds}
A proper upper bound for $\omega_{X/Y}^2$ implies the Effective Mordell conjecture. Note that $\tg \leq \frac{1}{4} \elg$ for a pm-graph $\ga$ containing bridges. Rumely showed that this inequality is sharp for a pm-graph $\ga$ with no cycles (i.e., the corresponding metrized graph has genus $0$). If $\ga$ is a bridgeless pm-graph, then $\tg \leq  \frac{1}{12} \elg$ by \cite[Corollary 5.8]{C2}, which is sharp for a circle graph as shown in \corref{cor tau for circle}. Moreover, $\tcg \leq 8 (\gc-1)^2 \tg$ for a pm-graph $\ga$ with genus $\gc \geq 2$ by \cite[Theorem 4.18]{C1}.
\end{remark}
\begin{remark}\label{rem valence vg}
Let $\ga$ be a graph with a vertex set $\#\vv{\ga}$. If we enlarge $\#\vv{\ga}$ by considering more valence $2$ points $p \in \ga$ assigned with $\bq (p)=0$ as vertices,
then we have the following observations:

The tau constant $\tg$ does not change, by its valence property
\cite[Remark 2.10]{C2}. The resistance function on $\ga$ does not change. The genus $\gc$ remains the same.
Moreover, $\va(p)-2+2 \bq (p)=0$ for each new vertex $p$. Therefore,
$\vg$ does not change. We call this property the \textit{valence property} of $\vg$.
\end{remark}
\begin{remark}\label{rem scale-independence vg}
For any given pm-graph $\ga$, if we multiply each edge length by a constant $t$, then the tau constant $\tg$ and
$r(x,y)$ for each $x$ and $y$ in $\ga$ change by a factor of $t$ \cite[Remark 2.15]{C2}.
This implies that $\vg = \varphi (\ga^N) \elg$ for any pm-graph $\ga$, where $\ga^N$ is the pm-graph obtained from
$\ga$ by dividing each edge length in $\ee{\ga}$ by $\elg$. We call this property the \textit{scale-independence} of $\vg$.
\end{remark}
\lemref{lem Zhang's second conj lower bound} shows how the first inequality in
\conjref{Conj Zhang's second} can be expressed in terms of invariants of bridgeless pm-graphs.
\begin{lemma}\label{lem Zhang's second conj lower bound}
The first inequality in \conjref{Conj Zhang's second} is equivalent to the following inequality:
\begin{equation*}
\begin{split}
\elg \leq 12 \tg +\frac{\tcg}{\gc -1}.
\end{split}
\end{equation*}
\end{lemma}
\begin{proof}
Both $\aga$ and $\ed$ are expressed in terms of $\tg$ and $\tcg$ in \propref{prop adm int of res}
and \propref{prop ed and tau}, respectively. After substituting these values into the first inequality in
\conjref{Conj Zhang's second} and doing some algebra, we obtain the result in the Lemma.
\end{proof}
\lemref{lem Zhang's second conj upper bound} shows that proving \conjref{conj Zhang varphi} is enough to show that the second equality in \conjref{Conj Zhang's second} holds.
\begin{lemma}\label{lem Zhang's second conj upper bound}
The second inequality in \conjref{Conj Zhang's second} is equivalent to the following inequality:
\begin{equation*}
\begin{split}
\vg \geq \frac{\cgc \elg}{4}.
\end{split}
\end{equation*}
\end{lemma}
\begin{proof}
The result follows by arguments similar to those in the proof of \lemref{lem Zhang's second conj lower bound}, using \thmref{thm main1}.
\end{proof}
\begin{theorem}\label{thm Zhang's second conj lower bound}
Let $\ga$ be a bridgeless pm-graph. Then we have
$$\frac{\gc-1}{\gc+1}(\ell(\ga)-4 \gc \cdot \aga) \leq \ed.$$
In particular, the first inequality in \conjref{Conj Zhang's second} holds.
\end{theorem}
\begin{proof}
Recall that $\ed \geq 0$ for any pm-graph $\ga$. Thus the inequality clearly holds when $\gc=1$.

Suppose $\gc \geq 2$. We first multiply both sides of the equality in \thmref{thm tau formula new} by $(\va(p)-2+2 \bq (p))$. Then we sum the resulting equality over all $p \in \vv{\ga}$, and use the fact that $\deg (K)=2 \gc -2$. In this way, we obtain
\begin{equation*}
\begin{split}
(2 \gc -2) \tg & = \frac{(2 \gc -2) \elg}{12}-\frac{1}{6}\sum_{p, \, q \in \vv{\ga}} (\va(q)-2)(\va(p)-2+2 \bq (p)) r(p,q)
\\ & \qquad +\frac{1}{3}\sum_{p \in \vv{\ga}} (\va(p)-2+2 \bq (p)) \sum_{e_i \in \, \ee{\ga}} \frac{\li}{\li+\ri}R_{c_i,p}.
\end{split}
\end{equation*}
Since $(\va(q)-2)=(\va(q)-2+2 \bq (q)-2 \bq (q))$,
\begin{equation*}
\begin{split}
(2 \gc -2) \tg & = \frac{(2 \gc -2) \elg}{12}-\frac{\tcg}{6}+\frac{1}{3}\sum_{p, \, q \in \vv{\ga}}(\va(p)-2+2 \bq (p)) \bq (q) r(p,q)
\\ & \qquad +\frac{1}{3}\sum_{p \in \vv{\ga}} (\va(p)-2+2 \bq (p)) \sum_{e_i \in \, \ee{\ga}} \frac{\li}{\li+\ri}R_{c_i,p}.
\end{split}
\end{equation*}
Equivalently,
\begin{equation}\label{eqn new form3}
\begin{split}
12 \tg +  & \frac{\tcg}{\gc -1}  = \elg + \frac{2}{\gc -1}\sum_{p, \, q \in \vv{\ga}}(\va(p)-2+2 \bq (p)) \bq (q) r(p,q)
\\ & \qquad +\frac{2}{\gc -1}\sum_{p \in \vv{\ga}} (\va(p)-2+2 \bq (p)) \sum_{e_i \in \, \ee{\ga}} \frac{\li}{\li+\ri}R_{c_i,p}.
\end{split}
\end{equation}
Since $K$ is effective and $\bq (p) \geq 0$ for each $p \in \vv{\ga}$, \eqnref{eqn new form3} implies that
\begin{equation}\label{eqn new form4}
\begin{split}
12 \tg + \frac{\tcg}{\gc -1} \geq \elg.
\end{split}
\end{equation}
Therefore, the inequality we wanted to show follows from \eqnref{eqn new form4} and \lemref{lem Zhang's second conj lower bound}.

Recall that an irreducible metrized graph is bridgeless. Hence, the first inequality in \conjref{Conj Zhang's second} holds.
\end{proof}
\begin{proposition}\label{prop main1}
Let $\ga$ be a bridgeless pm-graph. Then we have
\begin{equation*}
\begin{split}
\vg & = \frac{(2 \gc +1) \tg}{\gc}-\frac{\elg}{4 \gc}+ \frac{1}{2 \gc}\sum_{p, \, q \in \vv{\ga}}(\va(p)-2+2 \bq (p)) \bq (q) r(p,q)
\\ & \qquad +\frac{1}{2 \gc}\sum_{p \in \vv{\ga}} (\va(p)-2+2 \bq (p)) \sum_{e_i \in \, \ee{\ga}} \frac{\li}{\li+\ri}R_{c_i,p}.
\end{split}
\end{equation*}
In particular,
\begin{equation*}
\begin{split}
\vg \geq \frac{(2 \gc +1) \tg}{\gc}-\frac{\elg}{4 \gc}.
\end{split}
\end{equation*}
\end{proposition}
\begin{proof}By multiplying both sides of \eqnref{eqn new form3} by $\frac{\gc-1}{4 \gc}$ and using \thmref{thm main1}, we obtain the equality. Since $K$ is effective and $\bq$ is non-negative, the inequality follows.
\end{proof}
\lemref{lem Zhang's second conj lower bound}, the proof of \thmref{thm Zhang's second conj lower bound} and \propref{prop lambda} below clarify the relation between \conjref{conj Zhang lambda} and the first inequality in \conjref{Conj Zhang's second}.
\begin{proposition}\label{prop lambda}
Let $\ga$ be a bridgeless pm-graph. Then we have
\begin{equation*}
\begin{split}
\lag & = \frac{\gc}{8 \gc +4} \elg + \frac{1}{8 \gc +4} \sum_{p, \, q \in \vv{\ga}}(\va(p)-2+2 \bq (p)) \bq (q) r(p,q)
\\ & \qquad + \frac{1}{8 \gc +4} \sum_{p \in \vv{\ga}} (\va(p)-2+2 \bq (p)) \sum_{e_i \in \, \ee{\ga}} \frac{\li}{\li+\ri}R_{c_i,p}.
\end{split}
\end{equation*}
In particular,
\begin{equation*}
\begin{split}
\lag \geq \frac{\gc}{8 \gc +4} \elg.
\end{split}
\end{equation*}
\end{proposition}
\begin{proof}
 If we substitute the value of $\vg$ given in \propref{prop main1} and  the value of $\ed$ given in \propref{prop ed and tau} into the formula for $\lag$ given in \eqnref{eqn app3}, we obtain the following equality:
 \begin{equation}\label{eqn lambda1}
\begin{split}
\lag & = \frac{\gc -1}{2 \gc} \tg+ \frac{\tcg}{24 \gc}+\frac{4 \gc^2+\gc+1}{24 \gc (2 \gc +1)} \elg
\\ & \qquad +\frac{\gc-1}{12 \gc (2 \gc +1)}\sum_{p, \, q \in \vv{\ga}}(\va(p)-2+2 \bq (p)) \bq (q) r(p,q)
\\ & \qquad + \frac{\gc-1}{12 \gc (2 \gc +1)}\sum_{p \in \vv{\ga}} (\va(p)-2+2 \bq (p)) \sum_{e_i \in \, \ee{\ga}} \frac{\li}{\li+\ri}R_{c_i,p}.
\end{split}
\end{equation}
The equality in the proposition follows by multiplying \eqnref{eqn new form3} by $\frac{\gc -1}{24 \gc}$ and using \eqnref{eqn lambda1}.
Since $\ga$ is a pm-graph, the associated canonical divisor $K$ is effective. That is,
$\va(p)-2 + 2 \bq (p) \geq 0$ and $\bq (p) \geq 0$ for each $p \in \vv{\ga}$. Hence, the inequality in the proposition follows.
\end{proof}
\begin{corollary}\label{cor lambda for v 1 or 2}
Let $(\ga, \bq)$ be a bridgeless pm-graph with $\bq \equiv 0$. If $\ga$ has one or two vertices, then we have
\begin{equation*}
\begin{split}
\lag & = \frac{\gc}{8 \gc +4} \elg.
\end{split}
\end{equation*}
\end{corollary}
\begin{proof}
Let $\vv{\ga}= \{ p, q \}$. We have $R_{c_i,p}=0=R_{c_i,q}$, so the result follows by \propref{prop lambda}.
\end{proof}
\begin{remark}\label{rem lambda}
Whenever a pm-graph $(\ga,0)$ has two vertices say $p$ and $q$, we have $R_{c_i,p}=R_{c_i,q}=0$. Therefore, the lower bound for $\lag$ in \propref{prop lambda} is sharp. However, if $\ga$ has more than $3$ vertices, the bound can be far from sharp, as shown in \propref{prop fat necklace lambda}; and $\lag$ can be expressed as in \propref{prop lambda formula second}.
\end{remark}
In order to investigate the role of the case $\bq \equiv 0$ in finding lower bounds for $\vg$, $\lag$, $\ed$, and $\aga$, we make the following construction:

Let $(\ga,\bq)$ be a pm-graph of genus $\gc \geq 2$. If there is a vertex $p \in \vv{\ga}$ with $\bq (p) >0$, we attach $\bq (p)$ circles of length $\varepsilon>0$ to $\ga$ at the vertex $p$. By repeating this process for each such vertex, we obtain a new metrized graph, which we denote by $\ga_0$. By choosing $\bq_0 =0$ as the polarization on $\ga_0$, we have a pm-graph $(\ga_0,0)$. \figref{fig simplified} shows an example.
Note that $\vv{\ga_0}=\vv{\ga}$, $\ell(\ga_0)=\ell(\ga)+\varepsilon \sum_{p \in \vv{\ga}}\bq (p)$. Since $g(\ga_0)=g(\ga)+\sum_{p \in \vv{\ga}}\bq (p)$, $\gc (\ga_0)=\gc(\ga)$ by \eqnref{eqn genus}. Moreover,
 \begin{equation}\label{eqn simplified}
\begin{split}
\va_{\ga_0}(p) &=\va_{\ga}(p)+2 \bq(p), \quad \text{for each $p \in \vv{\ga}$, and}
\\ r_{\ga_0}(p,q)&=r_{\ga}(p,q), \quad \text{for each $p$ and $q$ in $\vv{\ga}$},
\end{split}
\end{equation}
where $\va_{\ga}(p)$ is valence of $p$ in $\ga$, and $r_{\ga}(x,y)$ is the resistance function on $\ga$. \eqnref{eqn simplified} implies
$\theta (\ga_0)=\tcg$. Moreover, $\tau(\ga_0)=\tg+\frac{\varepsilon}{12} \sum_{p \in \vv{\ga}}\bq (p)$ which can be seen by applying the additive property of the tau constant \cite[page 11]{C2}
and using the fact that $\tau(\beta)=\frac{\ell(\beta)}{12}$ for a circle graph $\beta$ (see \corref{cor tau for circle}). Thus, the following proposition follows from \thmref{thm main1}, \propref{prop ed and tau},
\propref{prop adm int of res}, and \corref{cor lambda interms of tau}.

\begin{proposition}\label{prop simplified}
Let $\ga$ be an arbitrary pm-graph with genus $\gc$, and let $\ga_0$ be as defined above. If $Q=\sum_{p \in \vv{\ga}}\bq (p)$, we have
\begin{align*}
\varphi(\ga_0) & = \vg + \varepsilon Q \frac{\gc -1}{6 \gc}, & \qquad \qquad
 \epsilon(\ga_0) & = \ed + \varepsilon Q \frac{\gc -1}{3 \gc},
\\ a(\ga_0)& = \aga+\varepsilon  Q \frac{2 \gc -1}{12 \gc^2} ,& \qquad \qquad
 \lambda(\ga_0)& = \lag+\varepsilon Q \frac{\gc}{8 \gc+4}.
\end{align*}
\end{proposition}
We obtain the following proposition by choosing an appropriate $\varepsilon$ as we construct $\ga_0$.
\begin{proposition}\label{prop simplified to epsilon}
Let $(\ga,\bq)$ be a given pm-graph with genus $\gc$. For any given $\varepsilon>0$, there exist a pm-graph $(\ga_0,0)$ of genus $\gc$ such that
$$\varphi(\ga_0) \leq \vg + \varepsilon, \qquad
 \epsilon(\ga_0) \leq \ed + \varepsilon, \qquad
a(\ga_0) \leq \aga+\varepsilon, \qquad
 \lambda(\ga_0) \leq \lag+\varepsilon.$$
\end{proposition}
\begin{remark}\label{rem simplified}
Since $\varepsilon$ in \propref{prop simplified to epsilon} can be taken arbitrarily small for a given pm-graph $\ga$, it will be enough to consider pm-graphs with polarization $\bq \equiv 0$ in order to give lower bounds for $\vg$, $\ed$, $\aga$, and $\lag$.
\end{remark}
Note that \propref{prop simplified to epsilon} should be compared with the article \cite[Lemma 5.14]{Fa} which gives a similar result for $\vg$.

\section{Simple polarized metrized graphs}\label{sec pm graph with zero part}

In this section, we will prove the second inequality in \conjref{Conj Zhang's second}, which is equivalent to a lower bound for $\vg$ by  \lemref{lem Zhang's second conj upper bound}. We will call a pm-graph  $(\ga, \bq)$ a \textit{simple pm-graph} if $\bq \equiv 0$. We will denote a simple pm-graph $(\ga,0)$ simply by $\ga$ when there is no danger of confusion. Note that $\gc =g$ (i.e., $\gc (\ga) =g(\ga)$) when $\bq (p)=0$ for each $p \in \vv{\ga}$ (see \eqnref{eqn genus}). To show that $\vg$ and $\lag$ are bounded by positive constants depending only on the genus and the length of the metrized graph, it will be enough to consider irreducible simple pm-graphs by \cite[4.4.2]{Zh2}, by our discussion about irreducible graphs in \secref{sec metrized graphs} and by \propref{prop simplified to epsilon}. Related results can be found in the article \cite[Lemmas 5.12, 5.14 and 5.15]{Fa}.
\begin{theorem}\label{thm varphi lower bounds}
Let $\ga$ be a bridgeless simple pm-graph with $\#(\vv{\ga})=v$ and $\#(\ee{\ga})=e$. Then we have
\begin{equation*}
\begin{split}
\vg & = \frac{(2 g +1) \tg}{g}-\frac{\elg}{4 g}+ \frac{1}{2 g}\sum_{p \in \vv{\ga}} (\va(p)-2) \sum_{e_i \in \, \ee{\ga}} \frac{\li}{\li+\ri}R_{c_i,p}.
\end{split}
\end{equation*}
Moreover,
\begin{align*}
&(i) \, \, \vg \geq \frac{(2 g +1) \tg}{g}-\frac{\elg}{4 g}, & \quad \quad \quad \quad \quad
& \quad \quad (ii) \, \, \vg \geq \frac{(g-1) \elg}{12 g (g+1)}, & \quad
\\ & (iii) \, \,   \vg \geq \frac{(4 e-5 v) \elg}{4 g (v+6)}, &  \quad \quad \quad \quad \quad
 & \quad \quad (iv) \, \, \vg \geq \frac{(g-1) \elg}{4 g (g+2)}, & \quad
\end{align*}
\begin{align*}
& (v) \, \, \vg \geq \frac{(2 g+1)g^2 v + 6 (2g+1 ) (v-1)^2-3e^2 v }{12g \cdot v \cdot e^2} \elg, \text{ if all edge lengths are equal}
\\ & \text{to each other and $\va(p) \geq 3$ for each $p \in \vv{\ga}$}.
\\ & (vi) \, \, \vg \geq \frac{2 g+1}{12g}\big(1-\frac{4}{\Lambda(\ga)} \big)^2 \elg+\frac{4(2 g+1)(\Lambda(\ga)-2)\elg}{g(v+6)\Lambda(\ga)^2}-\frac{\elg}{4g}, \text{ if } \Lambda(\ga) \geq 4.
\end{align*}
Therefore, the second inequality in \conjref{Conj Zhang's second} holds.
\end{theorem}
\begin{proof}
Since $\bq (p) =0$ for each $p \in \vv{\ga}$, the formula for $\vg$ in \propref{prop main1} reduces
to the formula given in the theorem.

Proof of $(i)$:
\\ It is given that $\ga$ is a pm-graph, so the associated canonical divisor $K$ is effective. That is,
in this case we have $\va(p)-2 \geq 0$ for each $p \in \vv{\ga}$. This gives $(i)$.

Proof of $(ii)$:
\\ We have $\tg \geq \frac{\elg}{6 (g+1)}$ by
\cite[Corollary 3.7]{C3}. Then the result follows from part $(i)$.

Proof of $(iii)$:
\\ We have $\tg \geq \frac{\elg}{2 (v +6)}$ by
\cite[Theorem 6.10 part (2)]{C3}. Then we have $\vg  \geq \frac{(2 g +1) \elg}{2 g (v+6)}-\frac{\elg}{4 g}$ by $(i)$ and we finish by using the fact that $g=e-v+1$.

Proof of $(iv)$:
Let $\vv{\ga}$ be a vertex set such that $\va (p) \geq 3$ if $p \in \ga$. Note that such a vertex set can be found by \remref{rem valence vg} for any bridgeless simple graph with genus at least two. Then by basic graph theory, $e \geq \frac{3 v}{2}$, where $e$ and $v$ are the number of edges and the number of vertices, respectively. Thus, $2(g-1) \geq v$. Then the result follows from $(iii)$.

Proof of $(v)$:
\\  When the edges have equal lengths and $\va(p) \geq 3$ for each $p \in \vv{\ga}$, we have $\tg \geq \big(\frac{1}{12}\big(\frac{g}{e}\big)^2+\frac{1}{2 v}\big(\frac{v-1}{e}\big)^2 \big) \elg$ by \cite[Theorem 6.11]{C3}
. Then the result follows from part $(i)$.

Proof of $(vi)$:
\\ When $\Lambda(\ga) \geq 4$, $\tg \geq \frac{1}{12} \big(1-\frac{4}{\Lambda(\ga)} \big)^2 \elg+\frac{4\elg(\Lambda(\ga)-2)}{(v+6)\Lambda(\ga)^2}$ by
\cite[Theorem 6.10 part (1)]{C3}. Then the result follows from part $(i)$.

Proof of the last part:
\\Using \lemref{lem Zhang's second conj upper bound},
proof of the second inequality in \conjref{Conj Zhang's second} follows from any of parts $(ii)$ and $(iv)$.
\end{proof}
Note that when $g \geq 4$, \thmref{thm main result2 on varphi} gives bounds to $\vg$ that are much stronger than the bounds given in \thmref{thm varphi lower bounds}.
\begin{remark}{\label{rem proof of Zhang second conj 1}}
The proof of \conjref{Conj Zhang's second} follows from
\thmref{thm Zhang's second conj lower bound} and \thmref{thm varphi lower bounds}.
\end{remark}
Suppose $\ga$ be a simple pm-graph with $\vv{\ga}= \{ p\}$ and $\# (\ee{\ga}) = e \geq 2$. We call such a graph a \textit{bouquet graph}. When $e=2$, $\ga$ is just a union of two circles along $p$.
\begin{proposition}\label{prop varphi for sum of circles}
Let $\ga$ be a simple bouquet graph. Then we have
\begin{equation*}
\begin{split}
\vg = \frac{g-1}{6 g} \elg.
\end{split}
\end{equation*}
\end{proposition}
\begin{proof}
The tau constant for a circle graph $\beta$ is $\frac{\ell (\beta)}{12}$ by \corref{cor tau for circle}. Then by the additivity of the tau constant \cite[page 11]{C2}
$\tg =\frac{\elg}{12}$. In this case, $R_{c_i,p}=0$ for any edge $e_i \in \ee{\ga}$. Since $\ga$ is a simple pm-graph, which means $\bq (p) =0$ for each $p \in \vv{\ga}$, the result follows from
the formula of $\vg$ in \propref{prop main1}.
\end{proof}
Suppose $\ga$ be a simple pm-graph with $\vv{\ga}= \{ p, \, q \}$ and $\# (\ee{\ga}) = e \geq 2$. We call such a graph a \textit{banana graph}. When $e=2$, $\ga$ is just a circle graph with two vertices. Note that
$g (\ga)=e-1$ for a banana graph with $e$ edges.
\begin{proposition}\label{prop varphi for banana}
Let $\ga$ be a simple banana graph with $e$ edges, and let the set of edge lengths be indexed by $\{L_1, L_2, \dots, L_e \}$. Then we have
\begin{equation*}
\begin{split}
\vg = \frac{g-1}{6 g} \elg -\frac{(g-1) (2g+1)}{6 g \sum_{i=1}^{e}\frac{1}{L_i}}.
\end{split}
\end{equation*}
In particular,
\begin{equation*}
\begin{split}
\vg \geq \frac{g (g-1)}{6 (g+1)^2} \elg.
\end{split}
\end{equation*}
\end{proposition}
\begin{proof}
The tau constant for banana graphs are calculated explicitly in the article \cite[Proposition 8.10]{C2}
. Namely, for a banana graph $\ga$, $\tg=\frac{\ell(\ga)}{12}-\frac{e-2}{6 \sum_{i=1}^{e}\frac{1}{L_i}}$. Since $\frac{\sum_{i=1}^e \li}{e}\geq \frac{e}{\sum_{i=1}^e\frac{1}{\li}}$ by Arithmetic-Harmonic Mean inequality, we have $\tg \geq \big( \frac{1}{12}-\frac{e-2}{6 e^2} \big) \elg$. In this case, $R_{c_i,q}=R_{c_i,p}=0$ for any edge $e_i \in \ee{\ga}$. Since $\ga$ is a simple pm-graph,
the result follows from
the formula of $\vg$ in \propref{prop main1}.
\end{proof}
When $\ga$ be a simple graph with two vertices, a formula for $\vg$ can be given explicitly by using \propref{prop varphi for banana} and the additivity of $\vg$.
\\ The following corollary shows that \thmref{thm varphi lower bounds} verifies or improves the previously known lower bounds to $\vg$ when $g \leq 4$.
\begin{corollary}\label{cor varphi bounds for small genus curves}
Let $\ga$ be a bridgeless simple pm-graph. Then $\vg \geq c(g) \cdot \elg$, where $c(g)$ can be taken as $\frac{1}{27}$ if $g=2$,
$\frac{1}{30}$ if $g=3$, $\frac{1}{32}$ if $g=4$, and $\frac{1}{35}$ if $g=5$.
%
\end{corollary}
\begin{proof}
We have $\vg \geq \frac{(g-1) \elg}{4 g (g+2)}$ by \thmref{thm varphi lower bounds} part $(iv)$. Evaluation of this lower bound for $g=3,4,5$ gives the desired lower bounds.

To obtain the given lower bound when $g=2$, we work with a vertex set such that $\va (p) \geq 3$ if $p \in \vv{\ga}$. Then, as discussed in the proof of part $(iv)$ in \thmref{thm varphi lower bounds}, we have $2(g-1) \geq v$. Therefore, $v \leq 2$ if $g=2$. Then either $\ga$ is the union of two circles along a point or a banana graph. Then the proof follows from \propref{prop varphi for sum of circles} and \propref{prop varphi for banana}.
\end{proof}
\begin{remark}{\label{rem varphi bounds for small genus curves}}
We recovered the known bound for $\vg$ when $g=2$; and improved the known bounds when $g=3,4$.
Note that the given lower bounds are sharp only for $g=2$.
\end{remark}
If any two distinct vertices of a graph $\ga$ are connected by one and only one edge and if there are no self loops, we call $\ga$ be a \textit{complete graph}.
For a complete graph on $v$ vertices, the valence of any vertex is $v-1$, and so by basic graph theory $e=\frac{v(v-1)}{2}$, and $g=\frac{(v-1)(v-2)}{2}$.
\begin{proposition}\label{prop varphi for complete graph}
Let $\ga$ be a simple complete pm-graph on $v$ vertices, where $v>3$. If all the edge lengths are equal, then we have
\begin{equation*}
\begin{split}
\vg = \frac{(v-2)(v-3)(v^2+6 v-6)}{6 v^3 (v-1)}\elg.
\end{split}
\end{equation*}
In particular,
$\frac{17}{288} \elg \leq \vg \leq \frac{9499}{54925} \elg$.
\end{proposition}
\begin{proof}
By \cite[Proposition 2.16]{C2}, $\tg =\Big( \frac{1}{12}\big(1-\frac{2}{v}\big)^2+\frac{2}{v^3}\Big)\ell (\ga)$. By arguments given in the proof of \cite[Proposition 2.16]{C2}, $r(p,q)=\frac{v-1}{e^2}$ for each pair of distinct vertices $p$, $q \in \vv{\ga}$. Therefore,
 $\tcg = (v-3)^2 \sum_{p, \, q \in \vv{\ga}} r(p,q)=\frac{4 (v-3)^2}{v}$. In addition, $e=\frac{v (v-1)}{2}$ and $g=\frac{(v-1)(v-2)}{2}$. Substituting these values into the formula for $\vg$ in \thmref{thm main1} gives the result. Then by elementary calculus we find that the minimum of $\vg$ is attained at $v=4$ and that its maximum is attained at $v=26$.
\end{proof}
By using the formula for $\lag$ given in \corref{cor lambda interms of tau} and following the arguments as in the proof of \propref{prop varphi for complete graph}, we can calculate $\lag$
for a simple complete pm-graph $\ga$.
\begin{proposition}\label{prop lambda for complete graph}
Let $\ga$ be a simple complete pm-graph on $v$ vertices, where $v>3$. If all the edge lengths are equal, then we have
\begin{equation*}
\begin{split}
\lag &= \frac{(v^3+v^2-12 v+18) (v-2)}{8 v^2 (v^2-3 v+3)} \elg
= \frac{g}{8 g+4}+\frac{(g+2) \sqrt{8 g+1}-7 g-2}{(2 g+1) \big(3+\sqrt{8 g+1}\big)^2}\elg.
\end{split}
\end{equation*}
In particular,
$\frac{25}{224} \elg \leq \lag \leq  \frac{499}{3650} \elg$, where the minimum is attained at $v=4$ and the maximum is attained at
$v=10$.
\end{proposition}

In view of \propref{prop lambda for complete graph}, there is not much room in improving the slope inequality given in (\ref{eqn slope inequality}).

We obtained a formula for $\tcg$ in \eqnref{eqn new form3}. The following Lemma gives another formula for $\tcg$. These two formulas will play important roles in giving an improved lower bound for $\vg$, which will give another proof of \conjref{conj Zhang varphi} and the second inequality in \conjref{Conj Zhang's second}.
\begin{lemma}\label{lem tcg second expansion}
Let $\ga$ be a bridgeless simple pm-graph with genus $g$ and $\# (\vv{\ga})=v$. Then we have
\begin{equation*}
\begin{split}
\tcg & = (2 g -2) \sum_{e_i \in \ee{\ga}} \frac{\li \ri}{\li +\ri}
+2 \sum_{q \in \vv{\ga}} (\va(q)-2)\sum_{e_i \in \ee{\ga}} \frac{R_{a_{i},q} R_{b_{i},q} + \ri R_{c_{i},q} }{\li +\ri}
\\ & \qquad +2 \sum_{q \in \vv{\ga}} (\va(q)-4)\sum_{e_i \in \ee{\ga}} \frac{\li R_{c_{i},q} }{\li +\ri} + 12 v \cdot \tg -v \cdot \elg.
\end{split}
\end{equation*}
\end{lemma}
\begin{proof}
First, we note the following equality for any given $q \in \vv{\ga}$.
\begin{equation}\label{eqn valence sum vertex-edge}
\begin{split}
\sum_{p \in \vv{\ga}} \va(p) r(p,q) & = \sum_{e_i \in \ee{\ga}} \big( r(\pp,q)+r(\qq,q) \big), \quad \text{$\pp$, $\qq$ are end points of $e_i$}.
\\ & = \sum_{e_i \in \ee{\ga}} \frac{\li \ri}{\li +\ri}+ 2 \sum_{e_i \in \ee{\ga}} \Big( \frac{R_{a_{i},q} R_{b_{i},q} }{\li +\ri} +R_{c_{i},q} \Big), \quad \text{by \eqnref{eqn2term0a}}.
\end{split}
\end{equation}
Summing the equality in \thmref{thm tau formula new} over all vertices gives
\begin{equation}\label{eqn middle sum over vertices}
\begin{split}
2 \sum_{p \, , q \in \vv{\ga}} (\va(q)-2) r(p,q)= v \cdot \elg - 12 v \cdot \tg + 4 \sum_{p \in \vv{\ga}} \sum_{e_i \in \ee{\ga}} \frac{\li R_{c_{i},p}}{\li + \ri}.
\end{split}
\end{equation}
On the other hand,
\begin{equation}\label{eqn tcg}
\begin{split}
\tcg & =\sum_{p \, , q \in \vv{\ga}} (\va(q)-2) (\va(p)-2) r(p,q), \quad \text{since $\ga$ is simple.}
\\ & =\sum_{q \in \vv{\ga}} (\va(q)-2) \sum_{p \in \vv{\ga}} \va(p) r(p,q) -2 \sum_{p \, , q \in \vv{\ga}} (\va(q)-2) r(p,q).
\end{split}
\end{equation}
Then, the result follows by substituting Equations (\ref{eqn valence sum vertex-edge}) and (\ref{eqn middle sum over vertices}) into \eqnref{eqn tcg} and using the fact that $\sum_{p \in \vv{\ga}} (\va(p)-2)=2g-2$.
\end{proof}
Next, we will combine \thmref{thm tau new formula with contraction} and \lemref{lem tcg second expansion} to obtain a new formula for $\tcg$.
\begin{proposition}\label{prop tcg third formula}
Let $\ga$ be a bridgeless simple pm-graph with genus $g$ and with at least three vertices. Then we have
\begin{equation*}
\begin{split}
\tcg & = -2 \elg + 24 \tg + (2 g -2) \sum_{e_i \in \ee{\ga}} \frac{\li \ri}{\li +\ri}
+ 4 \sum_{e_i \in \ee{\ga}} \frac{\ri}{\li +\ri} \sum_{e_j \in \ee{\oga_i}}  \frac{L_j \orcj}{L_j+\orj}
\\ & \qquad +2 \sum_{q \in \vv{\ga}} (\va(q)-4)\sum_{e_i \in \ee{\ga}} \frac{\li R_{c_{i},q} }{\li +\ri}.
\end{split}
\end{equation*}
\end{proposition}
\begin{proof}
Let $\# (\vv{\ga})=v$. We multiply both sides of the equality
given in \thmref{thm tau new formula with contraction} by $12 (v-2)$.
This will give us
\begin{equation*}
\begin{split}
2 \sum_{q \in \vv{\ga}}  (\va(q)-2) \sum_{e_i \in \ee{\ga}}  \frac{R_{a_{i},q} R_{b_{i},q}+\ri R_{c_i,q}}{\li+\ri} & =
 4 \sum_{e_i \in \ee{\ga}} \frac{\ri}{\li +\ri} \sum_{e_j \in \ee{\oga_i}}  \frac{L_j \orcj}{L_j+\orj}
\\ & \qquad \, +(v-2) \ell (\ga)-12 (v-2) \tg.
\end{split}
\end{equation*}
Then we substitute this into \lemref{lem tcg second expansion} to obtain the result.
\end{proof}
\begin{proposition}\label{prop varphi formula}
Let $\ga$ be a bridgeless simple pm-graph with genus $g$ and with at least three vertices. Then we have
\begin{equation*}
\begin{split}
\vg & = \frac{5 g +4}{g} \tg- \frac{g+2}{4 g} \elg + \frac{g -1}{2 g} \sum_{e_i \in \ee{\ga}} \frac{\li \ri}{\li +\ri}
+ \frac{1}{g} \sum_{e_i \in \ee{\ga}} \frac{\ri}{\li +\ri} \sum_{e_j \in \ee{\oga_i}}  \frac{L_j \orcj} {L_j+\orj}
\\ & \qquad +\frac{1}{2 g} \sum_{q \in \vv{\ga}} (\va(q)-4)\sum_{e_i \in \ee{\ga}} \frac{\li R_{c_{i},q} }{\li +\ri}.
\end{split}
\end{equation*}
In particular, if $\va(q) \geq 4$ for each $q \in \vv{\ga}$, then we have
\begin{equation*}
\begin{split}
\vg & \geq \frac{5 g +4}{g} \tg- \frac{g+2}{4 g} \elg + \frac{g -1}{2 g} \sum_{e_i \in \ee{\ga}} \frac{\li \ri}{\li +\ri}.
\end{split}
\end{equation*}
\end{proposition}
\begin{proof}
We obtain the equality by substituting the formula for $\tcg$ given in \propref{prop tcg third formula} into the formula for $\vg$ given in \thmref{thm main1}. Since all of $R_{c_{i},q}$, $\li$, $\ri$, $\orcj$, and $\orj$ are positive, the inequality follows from the equality if $\va(q) \geq 4$ for each $q \in \vv{\ga}$.
\end{proof}
Similarly, by using \propref{prop tcg third formula} and \corref{cor lambda interms of tau}, we obtain the following expression for $\lag$.
\begin{proposition}\label{prop lambda formula second2}
Let $\ga$ be a bridgeless simple pm-graph with genus $g$ and with at least three vertices. Then we have
\begin{equation*}
\begin{split}
\lag  & = \frac{3 g +3}{4 g+2} \tg + \frac{g-1}{16 g+8} \elg + \frac{g -1}{8 g+4} \sum_{e_i \in \ee{\ga}} \frac{\li \ri}{\li +\ri}
\\ & \quad+ \frac{1}{4 g+2} \sum_{e_i \in \ee{\ga}} \frac{\ri}{\li +\ri} \sum_{e_j \in \ee{\oga_i}}  \frac{L_j \orcj} {L_j+\orj}
 +\frac{1}{8 g+4} \sum_{q \in \vv{\ga}} (\va(q)-4)\sum_{e_i \in \ee{\ga}} \frac{\li R_{c_{i},q} }{\li +\ri}.
\end{split}
\end{equation*}
In particular, if $\va(q) \geq 4$ for each $q \in \vv{\ga}$, then we have
\begin{equation*}
\begin{split}
\lag & \geq \frac{3 g +3}{4 g+2} \tg + \frac{g-1}{16 g+8} \elg + \frac{g -1}{8 g+4} \sum_{e_i \in \ee{\ga}} \frac{\li \ri}{\li +\ri}.
\end{split}
\end{equation*}
\end{proposition}
The following is another formula for $\tcg$.
\begin{proposition}\label{prop tcg fourth formula}
Let $\ga$ be a bridgeless simple pm-graph with genus $g$ and with at least three vertices. Then we have
\begin{equation*}
\begin{split}
\tcg & = \frac{g-3}{2} \elg - 6 (g-3) \tg + (g -1) \sum_{e_i \in \ee{\ga}} \frac{\li \ri}{\li +\ri}
+ 2 \sum_{e_i \in \ee{\ga}} \frac{\ri}{\li +\ri} \sum_{e_j \in \ee{\oga_i}}  \frac{L_j \orcj}{L_j+\orj}
\\ & \qquad + 2 \sum_{q \in \vv{\ga}} (\va(q)-3)\sum_{e_i \in \ee{\ga}} \frac{\li R_{c_{i},q} }{\li +\ri}.
\end{split}
\end{equation*}
\end{proposition}
\begin{proof}
Since $\ga$ is simple, \eqnref{eqn new form3} can be expressed as follows:
\begin{equation}\label{eqn new form3 for simple}
\begin{split}
\tcg & = (g-1) \elg - 12 (g-1) \tg + 2\sum_{p \in \vv{\ga}} (\va(p)-2) \sum_{e_i \in \, \ee{\ga}} \frac{\li R_{c_i,p}}{\li+\ri}.
\end{split}
\end{equation}
Then the result is obtained by adding the formulas for $\tcg$ given in \propref{prop tcg third formula} and \eqnref{eqn new form3 for simple}.
\end{proof}
\begin{proposition}\label{prop varphi formula second}
Let $\ga$ be a bridgeless simple pm-graph with genus $g$ and at least three vertices. Then we have
\begin{equation*}
\begin{split}
\vg & = \frac{7 g +5}{2 g} \tg- \frac{g+3}{8 g} \elg + \frac{g -1}{4 g} \sum_{e_i \in \ee{\ga}} \frac{\li \ri}{\li +\ri}
+ \frac{1}{2 g} \sum_{e_i \in \ee{\ga}} \frac{\ri}{\li +\ri} \sum_{e_j \in \ee{\oga_i}}  \frac{L_j \orcj} {L_j+\orj}
\\ & \qquad +\frac{1}{2 g} \sum_{q \in \vv{\ga}} (\va(q)-3)\sum_{e_i \in \ee{\ga}} \frac{\li R_{c_{i},q} }{\li +\ri}.
\end{split}
\end{equation*}
In particular, if $\va(q) \geq 3$ for each $q \in \vv{\ga}$, then we have
\begin{equation*}
\begin{split}
\vg & \geq \frac{7 g +5}{2 g} \tg- \frac{g+3}{8 g} \elg + \frac{g -1}{4 g} \sum_{e_i \in \ee{\ga}} \frac{\li \ri}{\li +\ri}.
\end{split}
\end{equation*}
\end{proposition}
\begin{proof}
We obtain the equality by substituting the formula for $\tcg$ given in \propref{prop tcg fourth formula} into the formula for $\vg$ given in \thmref{thm main1}. Since all of $R_{c_{i},q}$, $\li$, $\ri$, $\orcj$, and $\orj$ are positive, the inequality follows from the equality if $\va(q) \geq 3$ for each $q \in \vv{\ga}$.
\end{proof}
Arguments similar to those in the proof of \propref{prop varphi formula second} can be used to derive another formula for $\lag$. Using \propref{prop tcg fourth formula} and \corref{cor lambda interms of tau}, we obtain
\begin{proposition}\label{prop lambda formula second}
Let $\ga$ be a bridgeless simple pm-graph with genus $g$ and with at least three vertices. Then we have
\begin{equation*}
\begin{split}
\lag  & = \frac{3 g +3}{8 g+4} \tg+ \frac{3 g-1}{16 (2 g+1)} \elg + \frac{g -1}{16 g+8} \sum_{e_i \in \ee{\ga}} \frac{\li \ri}{\li +\ri}
\\ & \quad+ \frac{1}{8 g+4} \sum_{e_i \in \ee{\ga}} \frac{\ri}{\li +\ri} \sum_{e_j \in \ee{\oga_i}}  \frac{L_j \orcj} {L_j+\orj}
 +\frac{1}{8 g+4} \sum_{q \in \vv{\ga}} (\va(q)-3)\sum_{e_i \in \ee{\ga}} \frac{\li R_{c_{i},q} }{\li +\ri}.
\end{split}
\end{equation*}
In particular, if $\va(q) \geq 3$ for each $q \in \vv{\ga}$, then we have
\begin{equation*}
\begin{split}
\lag & \geq \frac{3 g +3}{8 g+4} \tg + \frac{3 g-1}{16 (2 g+1)} \elg + \frac{g -1}{16 g+8} \sum_{e_i \in \ee{\ga}} \frac{\li \ri}{\li +\ri}.
\end{split}
\end{equation*}
\end{proposition}

Deriving new formulas for $\tg$ and $\tcg$ helped us obtain additional formulas for $\vg$ in \propref{prop tcg third formula} and \propref{prop tcg fourth formula}. We will show that these formulas lead to improved lower bounds for $\vg$. These new lower bounds are much more stronger than the ones given in \thmref{thm varphi lower bounds}. First, we will need to define the quantities below, which were used in first \cite[Section 3.9]{C1} and \cite{C3} to give lower bounds for $\tg$ and to establish a connection between $\tg$ and the edge connectivity $\Lambda(\ga)$.

%

In the rest of the paper, for any given pm-graph $\ga$ we will use the following notation $x(\ga)$ and $y(\ga)$, or simply $x$ and $y$ if there is no danger of confusion:
\begin{equation}\label{eqn edgecon1}
\begin{split}
x(\ga)&=\sum_{e_i \, \in
\ee{\ga}}\frac{\li^2\ri}{(\li+\ri)^2} +\frac{3}{4}\sum_{e_i \, \in
\ee{\ga}}\frac{\li \ri^2}{(\li+\ri)^2}
-\frac{3}{4}\sum_{e_i \, \in \ee{\ga}}\frac{\li
(R_{a_i,p}-R_{b_i,p})^2}{(\li+\ri)^2},
\\y(\ga)&=\frac{1}{4}\sum_{e_i \, \in \ee{\ga}}\frac{\li
\ri^2}{(\li+\ri)^2}+\frac{3}{4}\sum_{e_i \, \in \ee{\ga}}\frac{\li
(R_{a_i,p}-R_{b_i,p})^2}{(\li+\ri)^2}.
\end{split}
\end{equation}
If $\ga-e_i$ is not connected for an edge $e_i$, i.e. $\ri$ is infinite and $(R_{a_i,p}-R_{b_i,p})^2=\ri^2$, the
summands should be replaced with their corresponding limits as $\ri\longrightarrow \infty$.

By \eqnref{eqn edgecon1} and \propref{proptau}
\begin{equation}\label{eqn edgecon1 tau}
\begin{split}
\tg =\frac{\elg}{12}-\frac{x(\ga)}{6}+\frac{y(\ga)}{6},
\end{split}
\end{equation}
and it follows from \eqnref{eqn edgecon1} that
\begin{equation}\label{eqn edgecon1 x plus y}
\begin{split}
x(\ga)+y(\ga)=\sum_{e_i \in \ee{\ga}} \frac{\li \ri}{\li+\ri}.
\end{split}
\end{equation}
Next, we will express the bounds for $\vg$ found above in terms of $x(\ga)$ and $y(\ga)$.
\begin{proposition}\label{prop varphi lower bound for v(p) more than 3}
Let $\ga$ be a bridgeless simple pm-graph with genus $g$ and with at least three vertices.
Then we have the following two inequalities:
\\ If $\va (p) \geq 4$ for each $p \in \vv{\ga}$,
\begin{equation*}
\begin{split}
\vg & \geq \frac{g-1}{6 g} \elg-\frac{2 g +7}{6 g} x(\ga) + \frac{8 g +1}{6 g} y(\ga).
\end{split}
\end{equation*}
If $\va (p) \geq 3$ for each $p \in \vv{\ga}$,
\begin{equation*}
\begin{split}
\vg & \geq \frac{g-1}{6 g} \elg -\frac{g +2}{3 g} x(\ga) + \frac{5 g +1}{6 g} y(\ga).
\end{split}
\end{equation*}
\end{proposition}
\begin{proof}
By
\eqnref{eqn edgecon1 tau}
$\tg =\frac{\elg}{12}-\frac{x(\ga)}{6}+\frac{y(\ga)}{6}$, and by \eqnref{eqn edgecon1 x plus y} $x(\ga)+y(\ga)=\sum_{e_i \in \ee{\ga}} \frac{\li \ri}{\li +\ri}$. Substituting these into the inequalities given in Propositions \ref{prop varphi formula} and \ref{prop varphi formula second} gives the results.
\end{proof}
Let $\ga$ be a metrized graph such that $\va (p) =n \geq 2$ for each $p \in \vv{\ga}$. We call such a $\ga$ an \textit{n-regular} metrized graph, and we extend the definition to pm-graphs.

We consider a specific case in the following theorem.
\begin{theorem}\label{thm varphi for n-regular edge equal graphs}
Let $\ga$ be a bridgeless simple pm-graph with $\# (\vv{\ga})=v \geq 3$.
Suppose $\ga$ is n-regular and each edge in $\ee{\ga}$ has the same length.
If $n \geq 4$, we have
\begin{equation*}
\begin{split}
\frac{\vg}{\elg} \geq \frac{v^4 (n^2-4 n +10)(n-2)+4 v^3 (n^2-n-11)-v^2 (74 n -364)+12 v (5 n-43)+216}{6 n^2 v^3 (v (n-2)+2)}.
\end{split}
\end{equation*}
If $n=3$, we have
\begin{equation*}
\begin{split}
\vg & \geq \frac{4 v^4-8 v^3+91 v^2-222 v +144}{54 v^3 (v+2)} \elg.
\end{split}
\end{equation*}
\end{theorem}
\begin{proof}
Let $\# (\ee{\ga})=e$. By basic combinatorial graph theory, we have $e=\frac{n \cdot v}{2}$ for a $n$-regular metrized graph.
To make the proof simpler, we can assume that $\elg =1$ by using \remref{rem scale-independence vg}. Since the edge lengths are equal, $\li=\frac{1}{e}$ for each edge $e_i \in \ee{\ga}$. Then $x+y=\sum_{e_i \in \ee{\ga}}\frac{\li \ri}{\li+\ri}=\frac{1}{e}\sum_{e_i \in \ee{\ga}}\frac{ \ri}{\li+\ri}=\frac{(v-1)}{e}$ by \eqnref{eqn genusa}.
Then the result follows from
the inequalities in \propref{prop varphi lower bound for v(p) more than 3} and the inequality $y \geq \frac{v+6}{4 v} (x+y)^2$ (see the article \cite[Theorem 6.9 part (3)]{C3}).
\end{proof}
When $n=3$ and $v=4$, we have $\vg=\frac{17}{288}\elg$ for a simple pm-graph as in \thmref{thm varphi for n-regular edge equal graphs}. Therefore, the inequalities given in \thmref{thm varphi for n-regular edge equal graphs} are sharp.

By using relations between $x(\ga)$, $y(\ga)$, $\tg$ and $\Lambda(\ga)$, we will derive \thmref{thm main result1 on varphi} and \thmref{thm main result2 on varphi} which are the two main results on $\vg$:
\begin{figure}
\centering
\includegraphics[scale=0.6]{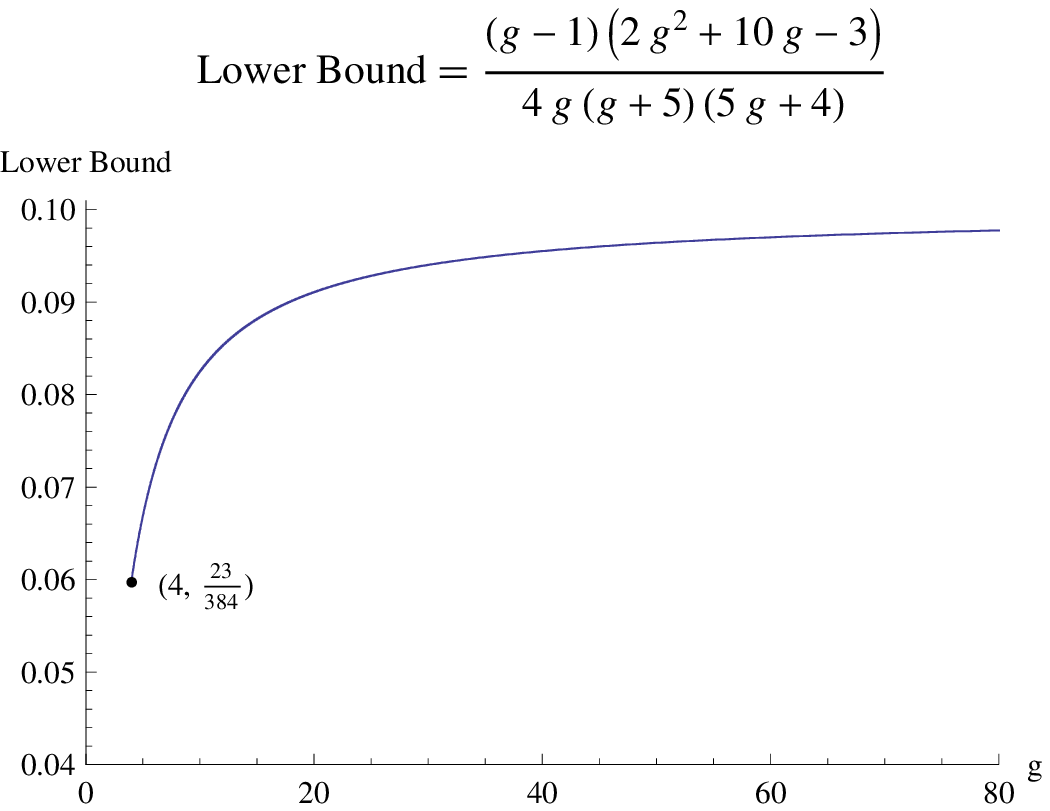} \caption{Lower bound to $\vg$ for $\ga$ as in \thmref{thm main result1 on varphi}.} \label{fig philowerbound1}
\end{figure}
\begin{theorem}\label{thm main result1 on varphi}
Let $\ga$ be a bridgeless simple pm-graph with genus $g$ and with $\# (\vv{\ga})=v \geq 3$.
Suppose $\va (p) \geq 4$ for each $p \in \vv{\ga}$. Then we have the following inequalities
\begin{equation*}
\begin{split}
\vg & \geq \frac{2 g^2 (v+10)-2 g (5v+2)-19v-16}{4 g (5 g+4) (v+6)} \elg.
\end{split}
\end{equation*}
In particular, we have (see \figref{fig philowerbound1}),
\begin{equation*}
\begin{split}
\vg & \geq \frac{(2 g^2 +10 g -3)(g-1) }{4 g (g+5) (5g+4)} \elg.
\end{split}
\end{equation*}
Moreover, if $\Lambda(\ga) \geq \frac{4 (5 g+4)}{2 g +7}$, the bounds above can be improved to
\begin{equation*}
\begin{split}
\vg & \geq \Big( \frac{(g-1) \Lambda(\ga)^2 (v+6)-4 v (2g+7)\Lambda(\ga) + 8v (5g+4)}{6g(v+6)\Lambda(\ga)^2}\Big)\elg
\end{split}
\end{equation*}
and
\begin{equation*}
\begin{split}
\vg & \geq \Big(\frac{g-1}{6 g} \big(1-\frac{4}{\Lambda(\ga)}\big)^2+
\frac{2 (g-1) (\Lambda(\ga)+2g-4)}{g(g+5)\Lambda(\ga)^2}\Big)\elg.
\end{split}
\end{equation*}
\end{theorem}
\begin{proof}
To make the proof simpler, we can assume that $\elg =1$ by using \remref{rem scale-independence vg}.
By \propref{prop varphi lower bound for v(p) more than 3} for our case, we have
\begin{equation*}
\begin{split}
\vg & \geq \frac{g-1}{6 g} \elg-\frac{2 g +7}{6 g} x + \frac{8 g +1}{6 g} y.
\end{split}
\end{equation*}
On the other hand, by \cite[Theorem 6.9]{C3}
, we have
$1\geq x+y$, $x \geq 0$, $y \geq 0$, $x \geq (\Lambda (\ga)-1) y$ and $y \geq \frac{v+6}{4 v} (x+y)^2$. Therefore, it will be enough to find $c=c(g,v)$ such that $\frac{g-1}{6 g}-\frac{2 g +7}{6 g} x + \frac{8 g +1}{6 g} y \geq c$ for any given fixed $g$ and $v$. We can choose $c$ such that $\frac{g-1}{6 g}-\frac{2 g +7}{6 g} x + \frac{8 g +1}{6 g} y = c$ is the tangent to the parabola $y = \frac{v+6}{4 v} (x+y)^2$. By elementary calculus, the tangency point will be given by $x_{0}=\frac{9 v (4 g^2+16 g+7) }{4 (v+6) (5 g+4)^2}$ and $y_{0}=\frac{v (2g+7)^2}{4 (v+6) (5 g+4)^2}$. These give $c=c(g,v)=\frac{2 g^2 (v+10)-2 g (5v+2)-19v-16}{4 g (5 g+4) (v+6)}$. This proves the first inequality in the theorem. An example is shown in \figref{fig philowerboundxy}.

By expressing $c(g,v)$ as a rational function in $v$, with rational coefficients involving $g$, we can rewrite $c(g,v)$ as
$\frac{2 g^2-10 g-19}{4 g (5g+4)}+\frac{8 g^2+56 g+98}{4g(5g+4)}\frac{1}{v+6}$. This is decreasing as a function of $v$ for each fixed $g$. Since $\va (p) \geq 4$ for each $p \in \vv{\ga}$, we have $g =e-v+1 \geq \frac{4 v}{2}-v+1=v+1$. Thus the second inequality follows by substituting $v=g-1$ into the first inequality.

Whenever $\Lambda(\ga) \geq \frac{4 (5 g+4)}{2 g +7}$, we can choose $c(g,v)$ such that
$\frac{g-1}{6 g}-\frac{2 g +7}{6 g} x + \frac{8 g +1}{6 g} y = c(g,v)$ passes through the point of intersection of  $x = (\Lambda (\ga)-1) y$ and $y= \frac{v+6}{4 v} (x+y)^2$. The final two results follow from this by elementary calculus.
\end{proof}
\begin{figure}
\centering
\includegraphics[scale=0.6]{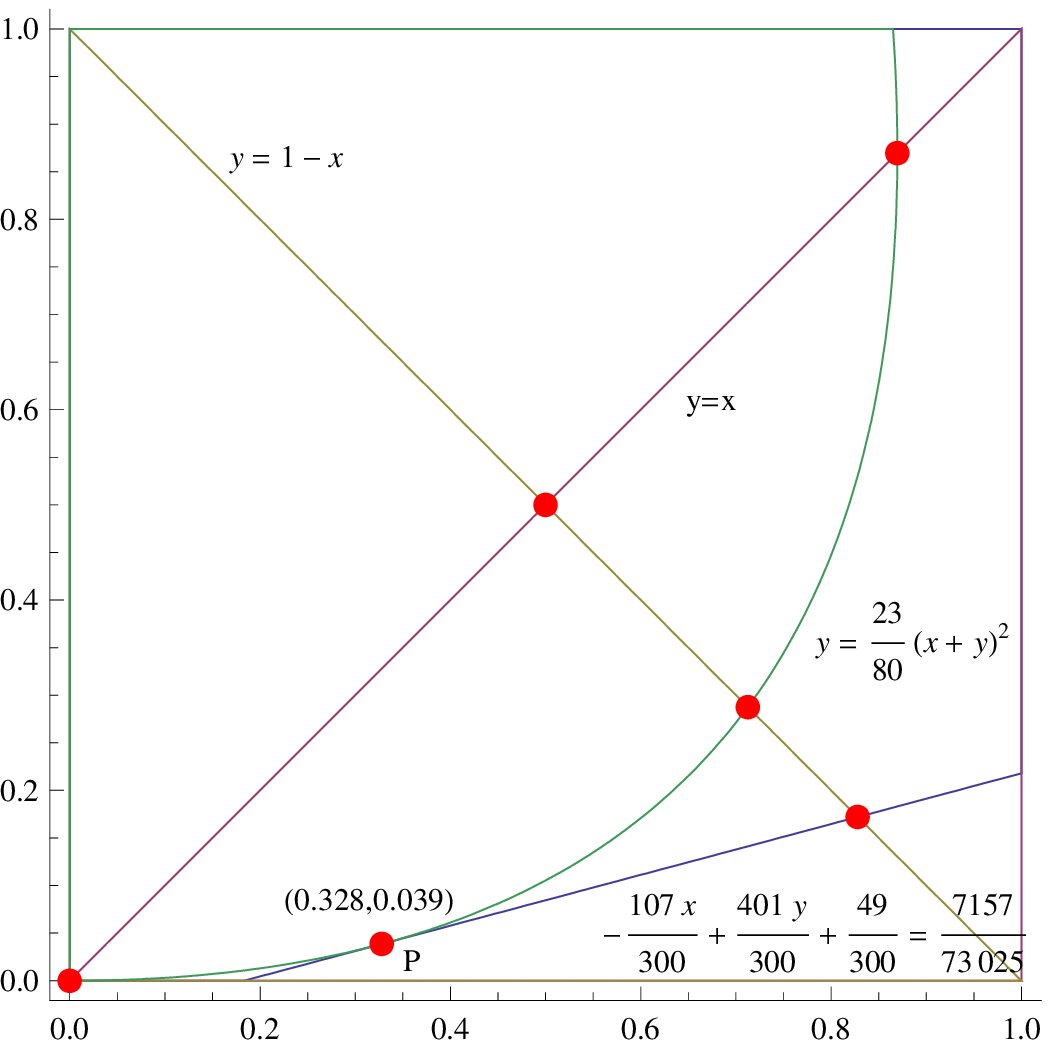} \caption{Relations between $x$ and $y$, and the lower bound to $\vg$, when $\ga$ is as in \thmref{thm main result1 on varphi} and $\Lambda (\ga)=2$, $g=50$, and $v=40$.} \label{fig philowerboundxy}
\end{figure}


\begin{figure}
\centering
\includegraphics[scale=0.6]{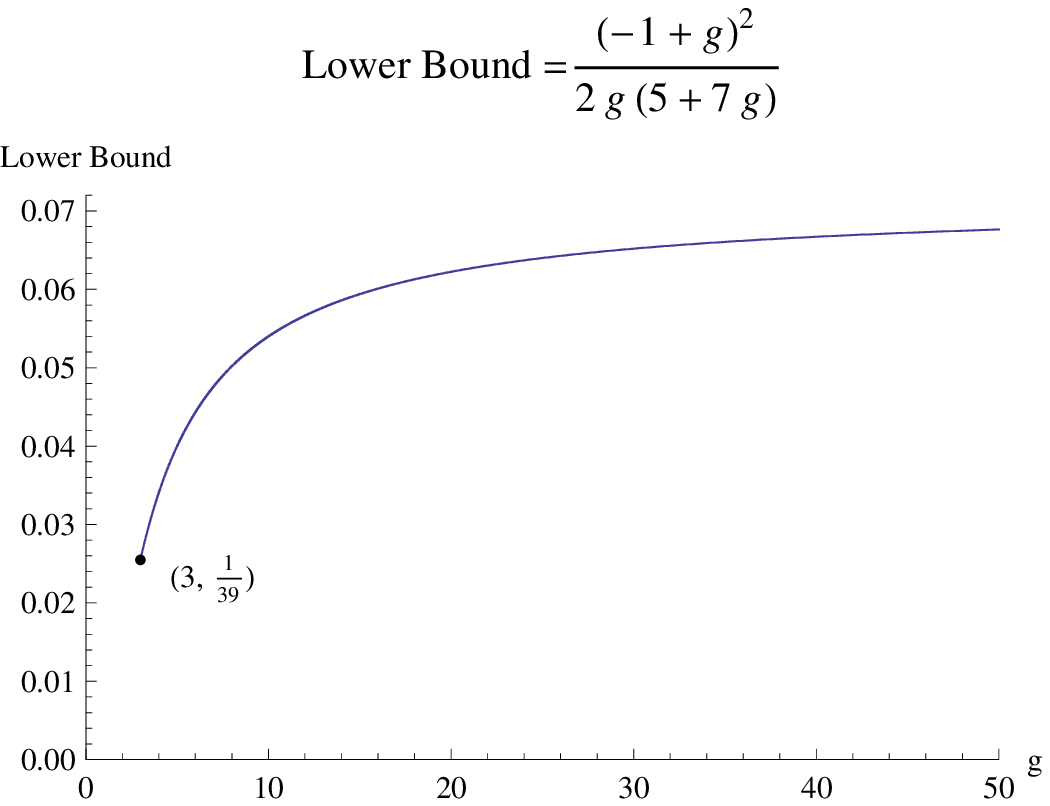} \caption{Lower bound to $\vg$ for $\ga$ as in \thmref{thm main result2 on varphi}.} \label{fig philowerboundforall}
\end{figure}
\begin{theorem}\label{thm main result2 on varphi}
Let $\ga$ be a bridgeless simple pm-graph with genus $g$ and with $\# (\vv{\ga})=v \geq 3$.
Then we have the following inequalities.
\\If $\va (p) \geq 3$ for each $p \in \vv{\ga}$, then we have
\begin{equation*}
\begin{split}
\vg & \geq \frac{g^2 (v +14) -2 g ( 3 v+2) -7 v-10}{2 g (7 g+5)(v+6)} \elg.
\end{split}
\end{equation*}
In particular, since $2(g-1) \geq v$ we have (see \figref{fig philowerboundforall}),
\begin{equation*}
\begin{split}
\vg & \geq \frac{(g-1)^2}{2 g (7 g+5)} \elg.
\end{split}
\end{equation*}
Moreover, if $\Lambda(\ga) \geq \frac{ 7 g+5}{g +2}$, the bounds above can be improved to
\begin{equation*}
\begin{split}
\vg & \geq \Big(\frac{(g-1) (v+6) \Lambda (\ga)^2-8v (g+2) \Lambda (\ga)+4 v (7g+5)}{6g(v+6)\Lambda(\ga)^2}\Big)\elg
\end{split}
\end{equation*}
and
\begin{equation*}
\begin{split}
\vg & \geq \Big(\frac{g-1}{6 g} \big(1-\frac{4}{\Lambda(\ga)}\big)^2+
\frac{2 (g-1)^2}{g(g+2)\Lambda(\ga)^2}\Big)\elg.
\end{split}
\end{equation*}
\end{theorem}
\begin{proof}
The proof follows from arguments similar to those in the proof of \thmref{thm main result1 on varphi}.
Note that we can use the second inequality in \propref{prop varphi lower bound for v(p) more than 3} for this time.
For the reader's convenience, we give some of the details.

We can assume that $\elg =1$ by using \remref{rem scale-independence vg}.
By \propref{prop varphi lower bound for v(p) more than 3} for this case, we have
\begin{equation*}
\begin{split}
\vg & \geq \frac{g-1}{6 g} \elg-\frac{g +2}{3 g} x + \frac{5 g +1}{6 g} y.
\end{split}
\end{equation*}
On the other hand, by \cite[Theorem 6.9]{C3}, we have
$1\geq x+y$, $x \geq 0$, $y \geq 0$, $x \geq (\Lambda (\ga)-1) y$ and $y \geq \frac{v+6}{4 v} (x+y)^2$. Therefore, it will be enough to find $c(g,v)$ such that $\frac{g-1}{6 g}-\frac{g +2}{6 g} x + \frac{5 g +1}{6 g} y \geq c(g,v)$ for any given fixed $g$ and $v$. We can choose $c=c(g,v)$ such that $\frac{g-1}{6 g}-\frac{g +2}{3 g} x + \frac{5 g +1}{6 g} y = c$ is the tangent to the parabola $y = \frac{v+6}{4 v} (x+y)^2$. By elementary calculus, the tangency point will be given by $x=\frac{12 v (2 g^2+5 g+2) }{(v+6) (7 g+5)^2}$ and $y=\frac{4 v (g+2)^2}{(v+6) (7 g+5)^2}$. These give $c(g,v)=\frac{g^2 (v+14)-2 g (3v+2)-7v-10}{2 g (7 g+5) (v+6)}$. This proves the first inequality in the theorem.

By expressing $c(g,v)$ as a rational function in $v$, with rational coefficients involving $g$, we can rewrite $c(g,v)$ as
$\frac{g^2-6 g-7}{2 g (7g+5)}+\frac{4 ( g^2+4 g+4)}{g(7g+5)}\frac{1}{v+6}$. This is decreasing as a function of $v$ for each fixed $g$. Since $\va (p) \geq 3$ for each $p \in \vv{\ga}$, $2 e =\sum_{p \in \vv{\ga}}\va (p) \geq 3 v$. Therefore, $g =e-v+1 \geq \frac{3 v}{2}-v+1=\frac{v}{2}+1$. Thus, the second inequality follows by substituting $v=2(g-1)$ into the first inequality.

Whenever $\Lambda(\ga) \geq \frac{(7 g+5)}{g +2}$, we can choose $c(g,v)$ such that
$\frac{g-1}{6 g}-\frac{g +2}{3 g} x + \frac{5 g +1}{6 g} y = c(g,v)$ passes through the point of intersection of  $x = (\Lambda (\ga)-1) y$ and $y= \frac{v+6}{4 v} (x+y)^2$. The final two results follow from this by elementary calculus.
\end{proof}
\begin{remark}{\label{rem minimal}}
For any given simple bridgeless pm-graph $\ga$ with genus $g \geq 2$ and a vertex set $\vv{\ga}$,
we can always find a non-empty vertex subset $\vv{\ga}':=\{ p \in \vv{\ga} \, | \, \va (p) \geq 3 \}$ by removing vertices of valence $2$ if there is any. By valence property of $\vg$ (\remref{rem valence vg}), this does not change $\vg$. We call $\vv{\ga}'$ be the \textit{minimal vertex set} of $\ga$. Note that $g \geq \frac{v}{2}+1$, where $\# (\vv{\ga}')=v$.
\end{remark}
\begin{theorem}\label{thm proof of Zhangs varphi conjecture}
Let $\ga$ be an irreducible pm-graph of genus $\gc \geq 2$. Then we have $\vg \geq t(\gc) \cdot \elg$, where $t(2)=\frac{1}{27}$, $t(3)=\frac{1}{30}$, and $t(\gc)=\frac{(\gc-1)^2}{2 \gc (7 \gc +5)}$ for $\gc \geq 4$.
\end{theorem}
\begin{proof}
This summarizes the best results we have shown above.

First, by considering \propref{prop simplified to epsilon} or equivalently \cite[Lemma 5.14]{Fa}, it will be enough to prove the desired lower bound inequalities for irreducible simple pm-graphs. Since every irreducible pm-graph is bridgeless, proving these lower bounds for bridgeless simple pm-graphs will implies that these lower bounds hold for irreducible pm-graphs.

Let $\ga$ be a bridgeless simple pm-graph. We can work with a minimal vertex set $\vv{\ga}$ by \remref{rem minimal}. Let $\# (\vv{\ga})=v$.

If $v=1$, then $\ga$ is a bouquet graph. By \propref{prop varphi for sum of circles}, we have $\vg = \frac{g-1}{6 g} \elg$, which is stronger than the desired results.

If $v=2$, then $\ga$ is a banana graph. By  \propref{prop varphi for banana}, we have $\vg \geq \frac{g (g-1)}{6 (g+1)^2} \elg$, which is stronger than what we wanted. When $g=2$, this gives $t(2)=\frac{1}{27}$.

If $v \geq 3$, we have $g \geq 3$ because we work with a minimal vertex set for $\ga$ (see \remref{rem minimal}). If $g=3$, we have $\vg \geq \frac{1}{30} \elg$ by \corref{cor varphi bounds for small genus curves}. If $g \geq 4$, we have $\vg \geq \frac{(\gc-1)^2}{2 \gc (7 \gc +5)} \elg$ by
\thmref{thm main result2 on varphi}. These finish the proof of the theorem.
\end{proof}
\begin{remark}\label{rem proof2}
By \lemref{lem Zhang's second conj upper bound}, \thmref{thm proof of Zhangs varphi conjecture} gives another proof of the second inequality in \conjref{Conj Zhang's second}, and improves the lower bounds given in \thmref{thm varphi lower bounds}. In \secref{sec genus 3}, we will consider $\gc=3$ case in more detail; and improve the lower bound $t(3)$. Namely, we can take $t(3)=\frac{892 - 11 \sqrt{79}}{14580} \approx 0.054473927$ by \thmref{thm g=3}.
\end{remark}
In \propref{prop lambda formula second}, we found an inequality for $\lag$ where $\ga$ is a simple pm-graph with more than $3$ vertices. Using \eqnref{eqn edgecon1 tau} and \eqnref{eqn edgecon1 x plus y}, this inequality can be stated in terms of $x(\ga)$ and $y(\ga)$
as follows:
\begin{equation}\label{lambda lower bound valence 3}
\begin{split}
\lag & \geq \frac{g}{8 g+4}+\frac{g y(\ga)-x(\ga)}{8 g+4}, \qquad \text{if $\va(p) \geq 3$ for each $p \in \vv{\ga}$.} 
\end{split}
\end{equation}
Similarly, it can be shown by using \propref{prop lambda formula second2} and Equations (\ref{eqn edgecon1}), (\ref{eqn edgecon1 tau}), and (\ref{eqn edgecon1 x plus y}) that
\begin{equation}\label{lambda lower bound valence 4}
\begin{split}
\lag &\geq \frac{g}{8 g+4}+\frac{g y(\ga)-x(\ga)}{4 g+2}, \qquad \text{if $\va(p) \geq 4$ for each $p \in \vv{\ga}$.}
\end{split}
\end{equation}

Note that by
\cite[Theorem 6.9 part (4)]{C3} we have  $g y \geq x \geq (\Lambda(\ga)-1) y > 0$, and recall that $y \geq \frac{v + 6}{4 v} (x + y)^2$, and $x+y < 1$ for any bridgeless metrized graph $\ga$ with $\elg =1$. In general, $y$ can be arbitrarily small. If all edge lengths are equal, then $x+y=\frac{v-1}{e}\elg$. Any ``proper'' improvement of these relations between $x$ and $y$ will result in improved lower bounds for $\tg$, $\lag$ and $\vg$. For various specific pm-graphs, it is easy to find the exact relation between $x$ and $y$.
%
%
\begin{figure}
\centering
\includegraphics[scale=0.5]{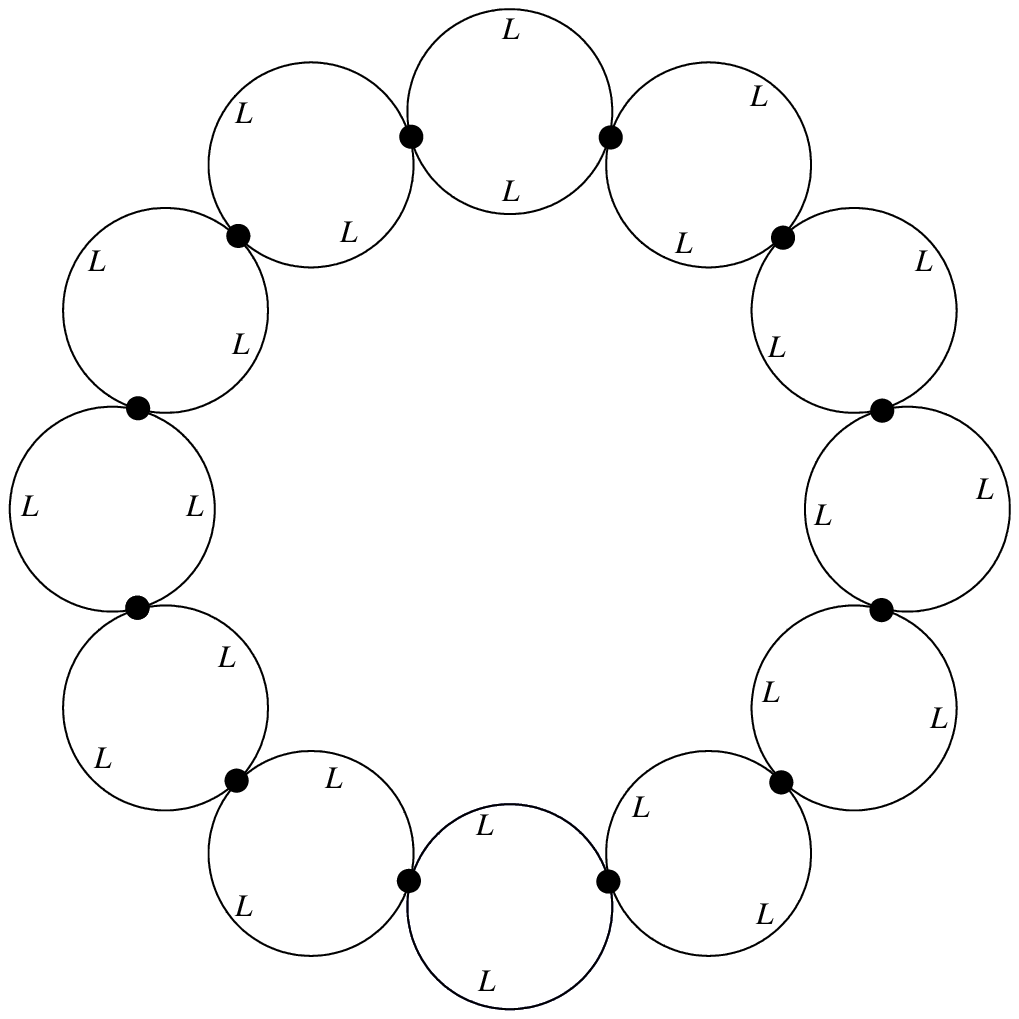} \caption{The graph $C_{12,2}$.} \label{fig circulardoublegraph3}
\end{figure}

We now give an example of a family of pm-graphs for which we can compute $\tg$, $\lag$ and $\vg$ by explicit formulas. This example shows how large $\lag$ and $\vg$ can be.
Let $C_{v,n}$ be the metrized graph obtained from the circle graph with $v$ vertices by replacing each of its edges by $n$ multiple edges of equal lengths so that the length of each edge in  $\ee{C_{v,n}}$ will be $\frac{\ell (C_{v,n})}{n \cdot v}$. \figref{fig circulardoublegraph3} illustrates an example. When $v=1$ and $v=2$, we gave formulas for $\lambda (C_{v,n})$ in \corref{cor lambda for v 1 or 2}, and for $\varphi (C_{v,n})$ in \propref{prop varphi for sum of circles} and \propref{prop varphi for banana}.
Suppose $v \geq 3$.
Let $\vv{C_{v,n}}=\{p_1,\, p_2, \dots, p_v \}$, and $r(x,y)$ be the resistance function on $C_{v,n}$. We have $\va(p) =2 n $ for each $p \in \vv{(C_{v,n}}$, so $e=n \cdot v$, and $g=(n-1) v +1$. Suppose all the edge lengths are equal. Then it follows from parallel and series circuit reductions that $r(p_v,p_i)=\frac{i (v-i)}{n^2 v^2} \ell (C_{v,n})$. Thus, $\sum_{i=1}^{v-1} r(p_v,p_i)= \frac{v^2 -1}{6 v n^2} \ell (C_{v,n})$. Moreover, by the symmetry of the graph, we have $\tcg = \sum_{i=1}^{v} \sum_{j=1}^{v} (\va(p_i)-2) (\va(p_j)-2) r(p_i,p_j)=\frac{2 (n-1)^2 (v^2 -1)}{3 n^2}\ell (C_{v,n})$. On the other hand, $\ta{C_{v,n}}=\big(\frac{(n-1)^2+1}{12 n^2}+\frac{n-1}{6 v n^2} \big) \elg$ by \cite[Example 3.9]{C2}. Hence, \corref{cor lambda interms of tau} and these results yield the following proposition:
\begin{proposition}\label{prop fat necklace lambda}
Let $C_{v,n}$ be a simple pm-graph with equal edge lengths, as defined before. For $v \geq 3$ we have
\begin{equation*}\label{}
\begin{split}
\lambda (C_{v,n}) = \frac{v^2 (n-1)^2 + 3 v (n-1) (n^2 -n +1) + 5 n^2 -4 n +2}{12 n^2 (2 (n-1)v+3)}\ell (C_{v,n}).
\end{split}
\end{equation*}
\end{proposition}
Using the above results along with \thmref{thm main1}, we obtain the following result.
\begin{proposition}\label{prop fat necklace varphi}
Let $C_{v,n}$ be a simple pm-graph with equal edge lengths, as defined before. For $v \geq 3$ we have
\begin{equation*}\label{}
\begin{split}
\varphi (C_{v,n}) = \frac{(n-1) ((v-2)^2+n v -1)}{6 n^2 v}\ell (C_{v,n}).
\end{split}
\end{equation*}
\end{proposition}
Note that for a fixed $v$, $\lim_{n\rightarrow \infty} \varphi (C_{v,n}) =\frac{1}{6}$, where length of the pm-graph $C_{v,n}$ is kept fixed as $1$.
\begin{corollary}\label{cor varphi necklace}
Let $C_{v,n}$ be a simple pm-graph with equal edge lengths, as defined before. For $v \geq 3$ we have
\begin{equation*}\label{}
\begin{split}
\varphi (C_{v,2}) & =\frac{(v-1)^2+2}{24 v} \ell(C_{v,2}) = \frac{(g-2)^2 +2}{24 (g-1)} \ell(C_{v,2}),
\\ \lambda (C_{v,2}) & =\frac{(v+2) (v+7)}{48 (2v+3)} \ell(C_{v,2}) = \frac{(g+6) (g+1)}{48 (2g+1)} \ell(C_{v,2}).
\end{split}
\end{equation*}
\end{corollary}
Note that \corref{cor varphi necklace} shows that $\vg$ and $\lag$ can be arbitrarily large for a graph of total length $1$.
%
\begin{figure}
\centering
\includegraphics[scale=0.55]{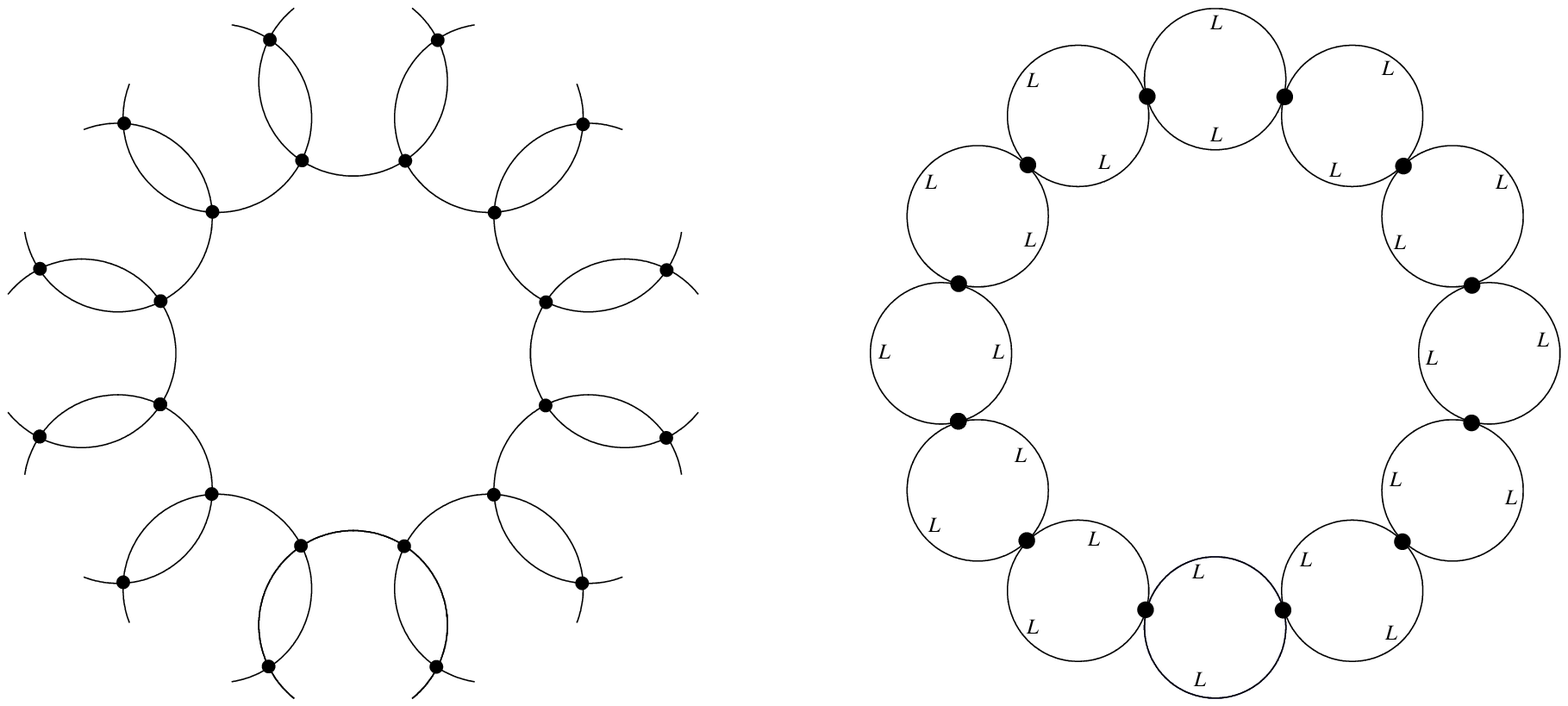} \caption{A curve with irreducible components and its dual graph.} \label{fig curvewd}
\end{figure}

Next, we will consider the simple complete pm-graphs with arbitrary length distributions,
unlike Propositions (\ref{prop varphi for complete graph}) and (\ref{prop lambda for complete graph}). Note that we have $g=\frac{(v-1)(v-2)}{2}$ for a simple complete pm-graph.
\begin{proposition}\label{prop varphi and lambda for complete graph}
Let $\ga$ be a simple complete pm-graph on $v$ vertices, where $v>3$. Then we have
$\tcg = 2 (v-3)^2 (x(\ga)+y(\ga))$,
$\vg = \frac{v (v-3) \elg}{6 (v-1)(v-2)}+\frac{v^2-21 v+48}{6 (v-1)(v-2)}x(\ga)+\frac{11 v^2-51 v+60}{6 (v-1)(v-2)}y(\ga)$ and
$\lag= \frac{(v - 1)^2 (v - 2)^2 - 2}{8((v - 1)^2 (v - 2)^2 - 1)} \elg
+ \frac{(v-2)(v-3) (v^2-7v+3)}{8((v - 1)^2 (v - 2)^2 - 1)}x(\ga)
+\frac{(v-3) (3 v^3- 15 v^2+ 23 v-6)}{8((v - 1)^2 (v - 2)^2 - 1)}y(\ga).$
\end{proposition}
\begin{proof}
We have $\va(p)=v-1$ for each $p \in \vv{\ga}$. On the other hand, $r(\pp,\qq)=\frac{\li \ri}{\li+\ri}$ if $\pp$ and $\qq$ are the end points
of the edge $e_i \in \ee{\ga}$. Therefore, $\tcg=2 (v-3)^2 \sum_{e_i \in \ee{\ga}}\frac{\li \ri}{\li+\ri}$. Then the formula for $\tcg$ follows from
\eqnref{eqn edgecon1 x plus y}.

The formula for $\vg$ follows from \eqnref{eqn edgecon1 tau}, the formula for $\tcg$ and \thmref{thm main1}.
Similarly, the formula for $\tcg$ follows from \eqnref{eqn edgecon1 tau}, the formula for $\tcg$ and \corref{cor lambda interms of tau}.
\end{proof}

\section{Calculations of pm-graph invariants using the discrete Laplacian}\label{sec calculations}

In this section, for a given pm-graph $\ga$ we express $\tg$ and $\tcg$ in terms of the corresponding discrete Laplacian. This leads to fast computer algorithms that can be used to compute $\vg$, $\ed$, $\lag$ (see the Mathematica codes given at the end).
Also, this gives another approach, as in \thmref{thm varphi inequality by lm}, to find lower bounds to these invariants.

The discrete Laplacian matrix of $\lm$ of a pm-graph $\ga$ is the same as the discrete Laplacian matrix of $\ga$ considered only with its metrized graph structure. This is no different than the construction of discrete Laplacian for a weighted graph. The details of its definition will be given below.

To have a well-defined discrete Laplacian matrix $\lm$ for a given pm-graph $\ga$, we first choose a vertex set $\vv{\ga}$ for $\ga$ in
such a way that there are no self-loops, and no multiple edges
connecting any two vertices. This can be done for any pm-graph $\ga$ by
enlarging the vertex set by inserting, whenever needed, additional valence two vertices
with $\bq$ value $0$. Note that enlarging the vertex set in this way does not change the value of $\vg$ (or $\ed$, $a(\ga)$, or $\lag$) by \remref{rem valence vg}. Such a vertex set $\vv{\ga}$, which can be used in the construction of a discrete Laplacian for $\ga$, will be called an \textit{optimal vertex set}. If two distinct vertices $p$ and $q$ are the end points of an edge, we call them \textit{adjacent vertices}.

\begin{definition}\label{def laplacian}
Let $\ga$ be a pm-graph with $e$ edges and with an optimal vertex set
$\vv{\ga}$ containing $v$ vertices. Fix an ordering of the vertices
in $\vv{\ga}$. Let $\{L_1, L_2, \cdots, L_e\}$ be a labeling of the
edge lengths. The $v \times v$ matrix $\am=(a_{pq})$ given by
\[
a_{pq}=\begin{cases} 0 & \quad \text{if $p = q$, or $p$ and $q$ are
not adjacent}.\\
\frac{1}{L_k} & \quad \text{if $p \not= q$, and $p$ and $q$ are
connected
by} \text{ an edge of length $L_k$}\\
\end{cases}.
\]
is called the \textit{adjacency matrix} of $\ga$. Let $\dm=\diag(d_{pp})$ be
the $v \times v$ diagonal matrix given by $d_{pp}=\sum_{s \in
\vv{\ga}}a_{ps}$. Then $\lm:=\dm-\am$ will be called the \textit{discrete
Laplacian matrix} of $\ga$. That is, $ \lm =(l_{pq})$ where
\[
l_{pq}=\begin{cases} 0 & \; \, \text{if $p \not= q$, and $p$ and $q$
are not adjacent}.\\
-\frac{1}{L_k} & \; \, \text{if $p \not= q$, and $p$ and $q$ are
connected by} \text{ an edge of length $L_k$}\\
-\sum_{s \in \vv{\ga}-\{p\}}l_{ps} & \; \, \text{if $p=q$}
\end{cases}.
\]
\end{definition}
We will denote the pseudo inverse of $\lm$ by $\plm$. For properties of $\lm$ and $\plm$, see the article \cite{C4} and the references given therein.
\begin{lemma} \cite{RB2} \cite{RB1} \label{lem disc2}
Let $\ga$ be a metrized graph with the discrete Laplacian $\lm$ and the
resistance function $r(x,y)$.
For the pseudo inverse $\plm$ we have
$$r(p,q)=l_{pp}^+-2l_{pq}^+ + l_{qq}^+, \quad \text{for any $p$, $q$ $\in \vv{\ga}$}.$$
\end{lemma}
Let $\ga$ be a pm-graph with an optimal vertex set $\vv{\ga}$. Using \lemref{lem disc2}, we can express $\tcg$ in terms of entries of $\plm$ as follows:
\begin{equation}\label{eqn tcg calculation}
\begin{split}
\tcg=\sum_{p, \, q \in \vv{\ga}} (\va(p)-2+2 \bq(p))(\va(q)-2+2 \bq(q)) (l_{pp}^+-2l_{pq}^+ + l_{qq}^+).
\end{split}
\end{equation}
In particular, suppose that $\va(p)-2+2 \bq(p)=k$ for each $p \in \vv{\ga}$. Since $2 \gc-2=\sum_{p \in \vv{\ga}} (\va(p)-2+2 \bq(p))$, $k=\frac{2 \gc-2}{v}$. On the other hand, we have $\sum_{p, \, q \in \vv{\ga}} (l_{pp}^+-2l_{pq}^+ + l_{qq}^+)=2 v \cdot  \text{trace}(\plm)$. This implies that
\begin{equation}\label{eqn tcg calculation regular}
\begin{split}
\tcg= \frac{8 (\gc-1)^2}{v} \text{trace}(\plm).
\end{split}
\end{equation}
\begin{theorem}\cite[Theorem 4.10]{C4}\label{thm disc2}
Let $\lm$ be the $v \times v$ discrete Laplacian matrix for a graph
$\ga$ with $v$ vertices. Let $\pp$ and $\qq$ be end points of edge $e_i$ for each $i
=1,2,\cdots,e$. Then
\begin{equation*}
\begin{split}
\tg
= -\frac{1}{12}\sum_{e_i \in \ee{\ga}} l_{\pp \qq}\big(\frac{1}{ l_{\pp \qq}}+l_{\pp \pp}^+-2l_{\pp
\qq}^+ +l_{\qq \qq}^+ \big)^2 -\frac{1}{4}\sum_{e_i \in \ee{\ga}} l_{p_i q_i} \big(l_{\pp \pp}^+ -  l_{\qq \qq}^+\big)^2 +\frac{1}{v}tr(\plm).
\end{split}
\end{equation*}
\end{theorem}
\begin{remark}\label{rem varphi algorithm}
For any given pm-graph $\ga$, we can choose an optimal vertex set $\vv{\ga}$. Then $\vg$ can be computed via a computer algorithm utilizing \thmref{thm main1}, \eqnref{eqn tcg calculation} and \thmref{thm disc2}.
\end{remark}
Since a pm-graph is connected, its discrete Laplacian matrix $\lm$ of size $v \times v$ has $v-1$ non-zero eigenvalues. Likewise,
its pseudo inverse $\plm$ has $v-1$ non-zero eigenvalues which are the reciprocal of the nonzero eigenvalues of $\lm$. Thus by applying Arithmetic-Harmonic Mean inequality, we have
\begin{equation}\label{eqn trace inequalities}
\begin{split}
\text{trace}(\plm)\geq \frac{(v-1)^2}{\text{trace}(\lm)}.
\end{split}
\end{equation}
For an $n$-regular metrized graph with equal edge lengths and with total length $1$, we have the following equality for the corresponding discrete Laplacian:
\begin{equation}\label{eqn trace formula for lm}
\begin{split}
\text{trace}(\lm)= v \cdot e \cdot n=\frac{v^2 \cdot n^2}{2}.
\end{split}
\end{equation}
and by \cite[Proof of Theorem 2.24]{C2}
\begin{equation}\label{eqn first term inequality}
\begin{split}
\sum_{e_i \in \ee{\ga}}\frac{\li^3}{(\li+\ri)^2} \geq \big(\frac{g}{e}\big)^2.
\end{split}
\end{equation}
We combine \eqnref{eqn trace formula for lm} and \eqnref{eqn trace inequalities} to obtain
\begin{equation}\label{eqn trace inequality for plm}
\begin{split}
\text{trace}(\plm) \geq \frac{2}{n^2} \big(\frac{v-1}{v}\big)^2.
\end{split}
\end{equation}
Therefore, by using \thmref{thm disc2}, the inequalities (\ref{eqn trace inequality for plm}) and (\ref{eqn first term inequality}), \thmref{thm main1}, and \remref{rem scale-independence vg}, we obtain the following inequality for $\vg$:
\begin{theorem}\label{thm varphi inequality by lm}
Let $(\ga,\bq)$ be a pm-graph such that $\bq(p)$ is constant for each $p \in \vv{\ga}$. Suppose that $\ga$ is a $n$-regular metrized graph with equal edge lengths
. Then we have
\begin{equation*}\label{}
\begin{split}
\vg \geq \frac{v^3 (n^2+2n-14) - v^2 (16 n -68) + 6 v (2n-15)+36}{6 n^2 v^3}\elg.
\end{split}
\end{equation*}
In particular, if $n=3$, we have
\begin{equation*}\label{}
\begin{split}
\vg \geq \frac{v^3+20v^2-54v+36}{54 v^3}\elg.
\end{split}
\end{equation*}
\end{theorem}
Note that when $v=4$ and $n=3$, $\vg =\frac{17}{288}$ for a simple graph as in \thmref{thm varphi inequality by lm}. Thus the inequalities given in \thmref{thm varphi inequality by lm} are sharp. Note also that the lower bounds given in \thmref{thm varphi inequality by lm} are independent of $\bq(p)$.

\section{Explicit formulas for genus 3 cubic simple pm-graphs}\label{sec genus 3}

In this section, we consider pm-graphs with genus $3$, and improve the lower bounds to $\vg$, $\lag$ by obtaining explicit formulas that are due to the techniques developed in the previous sections and the articles \cite{C2}, \cite{C3} and \cite{C4}.

\begin{figure}
\centering
\includegraphics[scale=0.7]{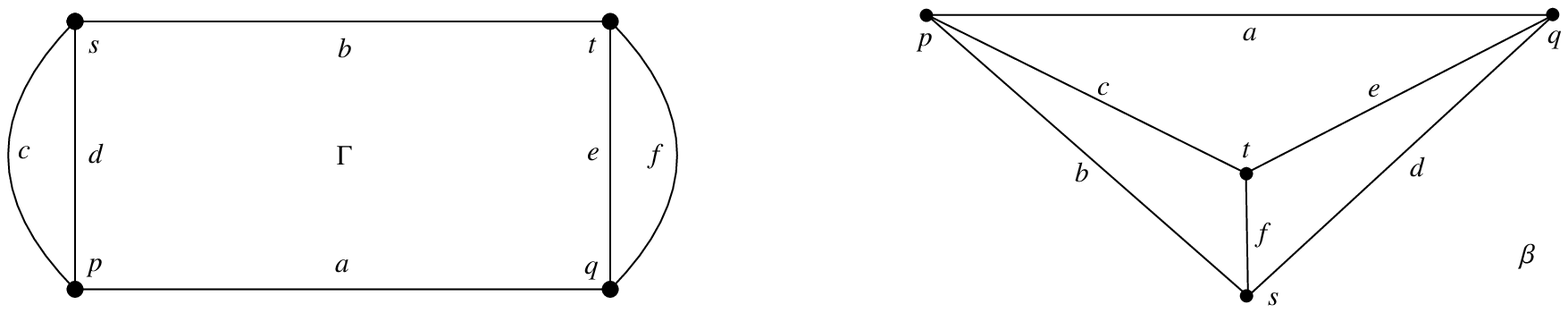} \caption{Cubic simple pm-graphs with genus $3$.} \label{fig g3curves}
\end{figure}

Recall that it is enough to find lower bounds for simple pm-graphs to find lower bounds to $\vg$, $\lag$
(see \remref{rem simplified}).
Moreover, Faber \cite[Lemma 5.15]{Fa} showed that it is enough to consider irreducible cubic simple pm-graphs to find lower bounds for $\vg$.
There are two types of bridgeless cubic simple pm-graphs of genus $3$ \cite[Figure 3]{Fa}. These are $\ga$ and $\beta$ as illustrated in \figref{fig g3curves}.

We will first consider $\ga$.
We have $\vv{\ga}=\{p, \, q, \, t, \, s \}$, $\ee{\ga}=\{ e_1, \, e_2, \, e_3, \, e_4, \, e_5, \, e_6 \}$ and the corresponding edge lengths
$\{a, \, b, \, c, \, d, \, e, \, f \}$ giving $\elg=a+b+c+d+e+f$.
Let $r_{\ga}(x,y)$ be the resistance function in $\ga$, and let $M=((c + d) (a + b) + c d) (e + f) + (c + d ) e f$.
Then by circuit reductions we have

$r_{\ga}(p,q)=a ((b (c + d) + c d) (e + f) + e f (c + d) )/M$,

$r_{\ga}(p,t)=(b c + b d + c d) (a e + a f + e f)/M$,

$r_{\ga}(p,s)=c d ((a + b + f) e + (a + b ) f)/M$,

$r_{\ga}(q,t)=((a + b + d) c + (a + b ) d) e f/M$,

$r_{\ga}(q,s)= (a (c + d) + c d) (b (e + f) + e f)/M$,

$r_{\ga}(s,t)=b (a (c + d) (e + f) + c d (e + f) + e f (c + d))/M$.

By \eqnref{eqn edgecon1 x plus y}, $x(\ga)+y(\ga)=r_{\ga}(p,q)+r_{\ga}(s,t)+2(r_{\ga}(p,s)+r_{\ga}(q,t))$. Thus,
\begin{equation}\label{eqn g3typeI x plus y}
\begin{split}
x(\ga)+y(\ga)=(3 (a + b) (c d (e + f) + (c + d) e f) + 2 a b (c + d) (e + f) + 4 c d e f)/M.
\end{split}
\end{equation}
Next, we give an explicit formula for the tau constant of $\ga$:
\begin{proposition}\label{prop g3typeI1}
Let $\ga$ be the simple pm-graph above. Then
\begin{equation*}
\begin{split}
\tg=\frac{\elg}{12}-((a + b) ((c + f) d e + c  f (d + e )) + 2 c d e f)/(6 M).
\end{split}
\end{equation*}
\end{proposition}
\begin{proof}
Recall that $e_1$ is an edge of $\ga$ with edge length $a$ as described above.
Applying \cite[Corollary 5.3]{C2} to $\ga$ with $e_1$, we obtain
$\tg=\ta{\ga-e_1}+\frac{a}{12}-\frac{R_1}{6}+\frac{A_{p,q,\ga-e_1}}{a+R_1}$. We have $\ta{\ga-e_1}=\frac{c+d+e+f}{12}+\frac{b}{4}$ by the additive property of the tau constant \cite[page 11]{C2}, \corref{cor tau for circle} and \cite[Corollary 2.22]{C2}. Moreover, $R_1=\frac{c d}{c+d}+b+\frac{e f}{e+f}$ by parallel and series circuit reductions, and $A_{p,q,\ga-e_1}=\frac{1}{6}\Big( \frac{c^2 d^2}{(c+d)^2}+\frac{e^2 f^2}{(e+f)^2}\Big)$ by \cite[Propositions 4.5, 4.6 and 8.9]{C2}. This gives the formula for $\tg$.
\end{proof}
%
%
\begin{proposition}\label{prop g3typeI2}
Let $\ga$ be the simple pm-graph above. Then
\begin{equation*}\label{}
\begin{split}
\tcg &=(6 (a + b) (c d (e + f) + (c + d) e f) + 8 a b (c + d) (e + f) +
 8 c d e f)/M,
\\ \vg &= \frac{\elg}{9}-(2 (a + b) (c d (e + f) + (c + d) e f) - 6 a b (c + d) (e + f) +
   7 c d e f)/(9 M).
\end{split}
\end{equation*}
\end{proposition}
\begin{proof}
Since $\tcg=2(r_{\ga}(p,q)+r_{\ga}(p,t)+r_{\ga}(p,s)+r_{\ga}(q,t)+r_{\ga}(q,s)+r_{\ga}(s,t))$, we have the formula for $\tcg$.
Since $\ga$ is simple, $\vg = \frac{(5 g -2) \tg}{g}+\frac{\tcg}{4 g}-\frac{\elg}{4}$
by \thmref{thm main1}. Then the result follows from the formula for $\tcg$ and \propref{prop g3typeI1}.
\end{proof}
%
%
%
%
%

%
By \eqnref{eqn edgecon1 tau}, \propref{prop g3typeI1} and \eqnref{eqn g3typeI x plus y} we have
\begin{equation}\label{eqn g3typeI x and y}
\begin{split}
x(\ga)&=(2 (a + b) (c d (e + f) + (c + d ) e f) + a b (c + d) (e + f) +
 3 c d e f)/M,
\\y(\ga) &=((a + b) (c d (e + f) + (c + d ) e f) + a b (c + d) (e + f) + c d e f)/M.
\end{split}
\end{equation}
\begin{proposition}\label{prop g3typeI3}
Let $\ga$ be the simple pm-graph above. Then we have the following sharp bound
$\vg > \frac{\elg}{16}$.
\end{proposition}
\begin{proof}
After doing the algebra we obtain that
$\vg = \frac{\elg}{16}+(11/2 (a + b) ( (c - d)^2 e + (c - d)^2 f + c (e - f)^2 +
d (e - f)^2) +
7 (d e (c - f)^2 + e f (d - c)^2 + c e (f - d)^2 + c f (e - d)^2 +
c d (e - f)^2 + d f (c - e)^2))/(144 M)+((14 (a^2 + b^2) + 220 a b) (c + d) (e + f) +
3 (a + b) ((c^2 + d^2 ) (e + f) + (c + d) ( e^2 + f^2)))/(288 M)$. This clearly gives that
$\vg > \frac{\elg}{16}$. If $a=b$ and $c=d=e=f$, then we have $\elg=2a+4c$, and
$\vg=\frac{\elg}{16}+\frac{a (62 a + 3 c)}{72 (2 a+c)}$. Since $\ga$ has genus $3$, we have $a>0$. Moreover,
$\vg$ approaches to $\frac{\elg}{16}$ as $a \rightarrow 0$.
\end{proof}
%
%

\begin{proposition}\label{prop g3typeI4}
Let $\ga$ be the simple pm-graph above. Then we have
$$\lag = \frac{3}{28}\elg+ (4 a b (c + d) (e + f) + (a + b) (c d (e + f) + (c + d) e f))/(28 M).$$
In particular, $\lag > \frac{3}{28}\elg$.
\end{proposition}
\begin{proof}
Since $\ga$ is simple, by \corref{cor lambda interms of tau} we have
$\lag = \frac{(3g -3) \tg}{4 g+2}+\frac{\tcg}{16 g+8}+\frac{(g +1)\elg}{16 g+8}$.
Hence, the result follows from \propref{prop g3typeI1} and the formula of $\tcg$ given in \propref{prop g3typeI2}.
\end{proof}
For $a=b$ and $c=d=e=f$, $\lag = \frac{3}{28}\elg+\frac{a}{14}$. Therefore, $\lag$ approaches to $\frac{3}{28}\elg$ as $a \rightarrow 0$.

Now, we will consider $\beta$.
We have $\vv{\beta}=\{p, \, q, \, s, \, t \}$, and the
edge lengths $\{a, \, b, \, c, \, d, \, e, \, f \}$ giving $\ell(\beta)=a+b+c+d+e+f$.
\[
\lm=\left[
\begin{array}{cccc}
 \frac{1}{a}+\frac{1}{b}+\frac{1}{c} & -\frac{1}{a} & -\frac{1}{b} & -\frac{1}{c}
 \medskip \\
 -\frac{1}{a} & \frac{1}{a}+\frac{1}{d}+\frac{1}{e} & -\frac{1}{d} & -\frac{1}{e}
 \medskip \\
 -\frac{1}{b} & -\frac{1}{d} & \frac{1}{b}+\frac{1}{d}+\frac{1}{f} & -\frac{1}{f}
 \medskip \\
 -\frac{1}{c} & -\frac{1}{e} & -\frac{1}{f} & \frac{1}{c}+\frac{1}{e}+\frac{1}{f}
\end{array}
\right],
\text{ and let }
\jm= \left[
\begin{array}{cccc}
 1 & 1 & 1 & 1 \\
 1 & 1 & 1 & 1 \\
 1 & 1 & 1 & 1 \\
 1 & 1 & 1 & 1
\end{array}
\right].
\]
Then, we compute the pseudo inverse $\plm$ of $\lm$ by using the formula $\plm=(\lm-\frac{1}{4} \jm)^{-1}+\frac{1}{4} \jm$. We have computed this by using Mathematica. Since the entries of $\plm$ are quite lengthy, we do not give the formula here. Let $N = a b d + a c d + b c d + a b e + a c e + b c e + b d e + c d e +
  a b f + a c f + b c f + a d f + c d f + a e f + b e f + d e f$.
We obtain the following formula for $\ta{\beta}$ by applying \thmref{thm disc2}:

$\ta{\beta}=\frac{\ell(\beta)}{12}-(c d e f + b (c d (2 e + f) + (c + d) e f) +
   a (c d (e + 2 f) + (c + d) e f) +
   a b (c (d + e + f) + d e + d f + 2 e f))/(6  N).$

We have $\theta (\beta) =8 \text{trace}(\plm)$ by \eqnref{eqn tcg calculation regular}, so

$\theta (\beta) =6 c d e f + b (8 c d e + 6 (c d + c e + d e) f) +
   a (c (6 d e + (8 d + 6 e ) f) + 6 d e f +
      b (6 c (d + e + f) + 6 d e + (6 d + 8 e ) f )))/ N.$

We have $\theta (\beta) = 2(x(\beta)+y(\beta))$ by \eqnref{eqn tcg calculation regular}, and
$\ta{\beta}=\frac{\ell(\beta)}{12}-\frac{x(\beta)}{6}+\frac{y(\beta)}{6}$ by \eqnref{eqn edgecon1 tau}.
These give

$y(\beta)=(a b c d + a b c e + a b d e + a c d e + b c d e + a b c f + a b d f +
   a c d f + b c d f + a b e f + a c e f + b c e f + a d e f +
  b d e f + c d e f)/N$, and
$x(\beta)=2 y(\beta)+(b c d e + a c d f + a b e f)/N$.

\begin{figure}
\centering
\includegraphics[scale=0.65]{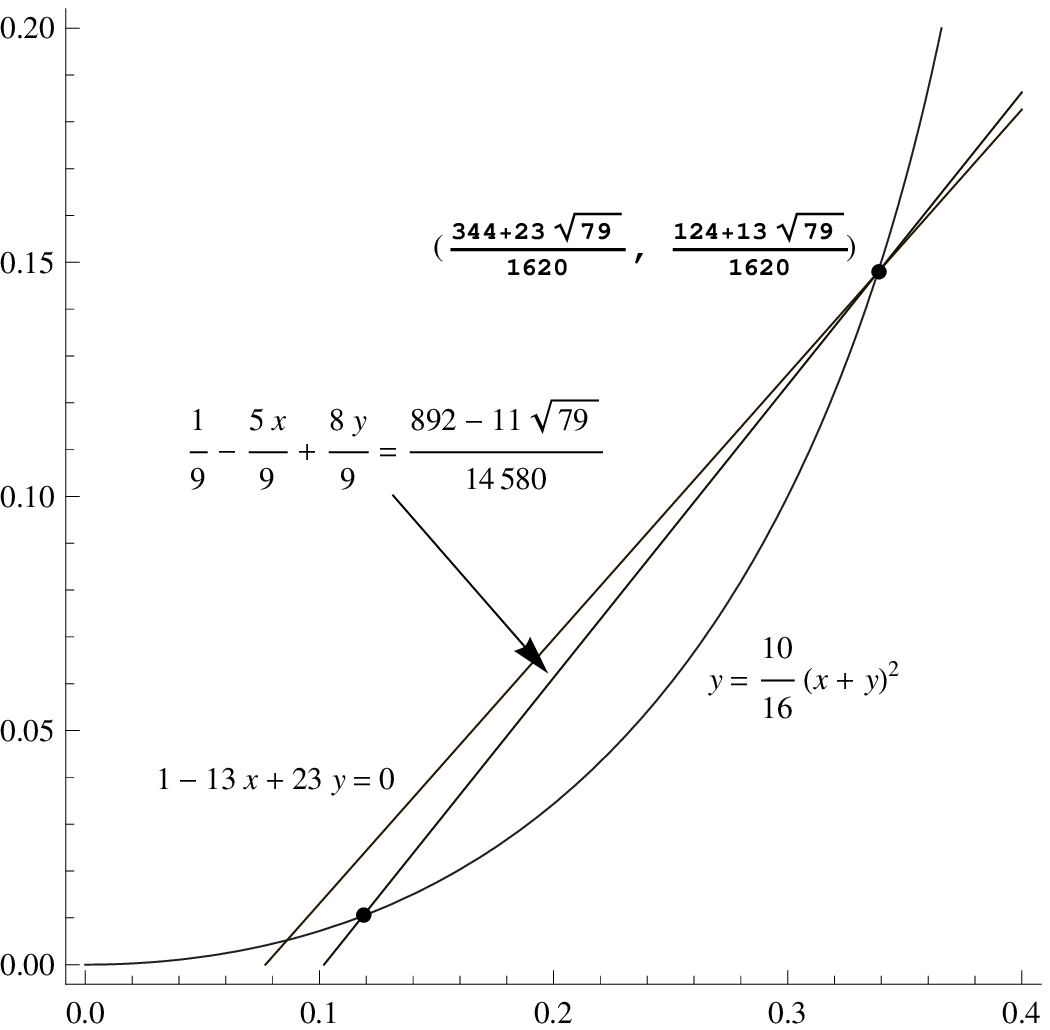} \caption{Giving bounds to $\varphi(\beta)$.} \label{fig g3graphbound}
\end{figure}

By \thmref{thm main1} and above formulas for $\ta{\beta}$ and $\theta(\beta)$, we have
$\varphi(\beta)= \frac{\ell(\beta)}{9}-\frac{5}{9} x(\beta)+\frac{8}{9}y(\beta)$.
Note that this shows that the second lower bound given in \propref{prop varphi lower bound for v(p) more than 3} is sharp.
\begin{proposition}\label{prop g3 beta lambda}
Let $\beta$ be the simple pm-graph above. Then
$\lambda(\beta) = \frac{3}{28}\ell(\beta)+ (c d e f + b (d e f + c (d f + e f)) +
 a (d e f + c (d e + e f) + b (d e + d f + c (d + e + f))))/(28 N)$.
In particular, $\lambda(\beta) > \frac{3}{28}\ell(\beta)$.
\end{proposition}
\begin{proof}
Since $\beta$ is simple, by \corref{cor lambda interms of tau} we have
$\lambda(\beta) = \frac{(3g -3) \ta{\beta}}{4 g+2}+\frac{\theta (\beta)}{16 g+8}+\frac{(g +1)\ell(\beta)}{16 g+8}$.
Then the result follows from the above formulas for $\ta{\beta}$ and $\theta (\beta)$.
\end{proof}
For $a=f=\frac{\ell(\beta)}{2}-k$ and $b=c=d=e=\frac{k}{2}$, $\lambda(\beta) = \frac{3}{28}\ell(\beta)+\frac{k(\ell(\beta)-2 k)}{56 (\ell(\beta)-k)^2}\ell(\beta)$. Therefore,
$\lambda(\beta)$ approaches to $\frac{3}{28}\ell(\beta)$ as $k \rightarrow 0$.
\begin{proposition}\label{prop g=3 beta}
Let $\beta$ be as before. Then
$13 x(\beta) \geq \ell(\beta)+23 y(\beta).$
\end{proposition}
\begin{proof}
Let $h$ be the function given by $h(a,b,c,d)=a^2 b c + a b^2 c + a b c^2 + a^2 b d + a b^2 d + a^2 c d -
 12 a b c d + b^2 c d + a c^2 d + b c^2 d + a b d^2 + a c d^2 +
 b c d^2$. We have Arithmetic-Harmonic mean inequality for any given $a>0$, $b>0$, $c>0$ and $d>0$. Namely,
 $a+b+c+d \geq \frac{16}{\frac{1}{a}+\frac{1}{b}+\frac{1}{c}+\frac{1}{d}}$ with equality iff $a=b=c=d$. Note that this is equivalent to $h(a,b,c,d) \geq 0$. By using the above formulas of $x(\beta)$ and $y(\beta)$, we obtain
 $\ell(\ga)-13 x(\ga)+23 y(\ga)=(h(b, c, d, e) + h(a, c, d, f) + h(a, b, e, f)+a^2 b d + a b^2 d + a b d^2 + a^2 c e + a c^2 e + a c e^2 + b^2 c f +
 b c^2 f + d^2 e f + d e^2 f + b c f^2 + d e f^2)/N$. This gives the result.
\end{proof}
\begin{proposition}\label{prop g=3 beta2}
Let $\beta$ be as before. Then
$\varphi(\beta) \geq \frac{892 - 11 \sqrt{79}}{14580} \ell(\beta) \approx 0.054473927 \ell(\beta)$.
\end{proposition}
\begin{proof}
We assume that $\ell(\beta)=1$. We have $y > \frac{10}{16} (x + y)^2$ by \cite[Theorem 6.9 part (3)]{C3}, and $13 x \geq \ell(\beta)+23 y$ by \propref{prop g=3 beta}. Note that the parabola  $y = \frac{10}{16} (x + y)^2$ and the line
$13 x = \ell(\beta)+23 y$ have intersection point with coordinates $x_{0}=\frac{344 + 23 \sqrt{79}}{1620}$ and
$y_{0}=\frac{124 + 13 \sqrt{79}}{1620}$ as also illustrated in \figref{fig g3graphbound}. These inequalities imply that $\frac{\ell(\beta)}{9}-\frac{5}{9} x +\frac{8}{9}y \geq \frac{892 - 11 \sqrt{79}}{14580}$. That is, $\varphi(\beta) \geq \frac{892 - 11 \sqrt{79}}{14580}$. Then the result follows from \remref{rem scale-independence vg}.
\end{proof}
When $a=b=c=d=e=f=\frac{1}{6}$, $\varphi(\beta)=\frac{17}{288}$.
There is considerable computational evidence, based on the numerical optimization methods in Mathematica, that the lower bound in \propref{prop g=3 beta2} can be taken as $\frac{17}{288} \ell(\beta)$, which is also conjectured by Faber \cite[Remark 5.1]{Fa}.

The lower bounds to $\vg$ were given in \thmref{thm proof of Zhangs varphi conjecture}. When $\gc =3$, the previous lower bound can be improved as follows:
\begin{theorem}\label{thm g=3}
Let $\ga$ be an irreducible pm-graph with genus $\gc=3$. Then we have $\vg \geq \frac{892 - 11 \sqrt{79}}{14580} \ell(\beta) \approx 0.054473927 \ell(\beta)$.
\end{theorem}
\begin{proof}
First, by considering \propref{prop simplified to epsilon} or equivalently \cite[Lemma 5.14]{Fa}, it will be enough to prove the given lower bound inequalities for any given irreducible simple pm-graphs. Since every irreducible pm-graph is bridgeless, proving these lower bounds for bridgeless simple pm-graphs will implies that these lower bounds hold for irreducible pm-graphs. Moreover,
by using \cite[Lemma 5.15]{Fa}, it will be enough to consider bridgeless simple cubic pm-graphs of genus $3$.
Recall that there are only two types of bridgeless cubic pm-graph of genus $g=3$.
Hence, the result follows from Propositions (\ref{prop g=3 beta2}) and (\ref{prop g3typeI3}).
\end{proof}
One can use the techniques that we have developed to compute $\vg$ and $\lag$ explicitly for other graphs, similarly to what we did in this section. However, the formulas become too large to include as $g$ gets bigger.








\tiny{
\begin{verbatim}

(* Mathematica Package Metrized Graph Invariants, i.e, MGI` in short. *)
(*                     written by Zubeyir Cinkir                      *)

(* All functions are Compatible with Mathematica Version 7, and all functions except for TotalLength are
compatible with Mathematica Version 6.*)

BeginPackage["MGI`"]

Valence::usage = "Valence[A] is a pure function that returns the valence of a given vertex of the metrized
graph \[CapitalGamma] with discrete Laplacian A. Therefore, Valence[A][p] is the valence of the vertex p
corresponding to the p-th row of A."

Genus::usage = "Genus[A] gives the GenusM[A] gives genus of the metrized graph \[CapitalGamma] corresponding
to the discerete Laplacian matrix A. This is the same as the 1st Betti number of \[CapitalGamma], which is
(number of edges in \[CapitalGamma]) - (number of vertices in \[CapitalGamma]) + 1."

CompleteGraphM::usage = "CompleteGraphM[n] gives the discrete Laplacian matrix of the complete graph on n
vertices with equal edge lengths such that the total length of the graph is 1."

TotalLength::usage = "TotalLength[A] gives the total length of the graph \[CapitalGamma] with discrete
Laplacian A. The total length of a graph is the sum of its edge lengths."

ResistanceMatrix::usage = "ResistanceM[A] gives the resistance matrix of the graph \[CapitalGamma] with the
discrete Laplacian matrix A. ResistanceM[[p,q]] is the effective resistance between the vertices p and q,
which correspond to the p-th and q-th rows of A. In this case, \[CapitalGamma] is considered as the
electric circuit such that the resistances along the edges are given by the edge lengths, and that unit
current enters at p and leaves at q."

AdjacentToLaplacian::usage = "AdjacentToLaplacian[A] AdjacentToLaplacian[A] constructs the discrete
Laplacian matrix of the graph with the Adjacency matrix A."

ThetaConstant::usage = "ThetaConstant[A,Q] computes
\[Theta](\[CapitalGamma])= \!\(\*UnderscriptBox[\" \[Sum] \", RowBox[{\" p \", \",\", \" \", RowBox[{\" q \"
, \" \", \" vertices \", \" \", \" in \", \" \", \[CapitalGamma] }]}]]\)(Valence[A][p] - 2 + 2 Q[[p]]) (
Valence[A][q]- 2 + 2 Q[[q]]) r (p, q), where r(x,y) is the resistance function on the graph \[CapitalGamma]
with discrete Laplacian A. Here Q is either 0 or a list of non-negative integers, in which case the size of
Q is equal to the number of rows in A."

TauConstant::usage = "TauConstant[A] gives the tau constant of the metrized graph \[CapitalGamma] with the
discrete Laplacian matrix A."

PhiConstant::usage = "PhiConstant[A,Q] computes \[CurlyPhi] (\[CapitalGamma]), where \[CapitalGamma] is the
metrized graph with the discrete Laplacian A. Here Q is either 0 or a list of non-negative integers, in
which case the size of Q is equal to the number of rows in A."

Begin["`Private`"]
(* Implementation of the package *)

(*Note that PseudoInverse[A] == Inverse[A - 1/k] + 1/k == Inverse[A +1/k] - 1/k, where k = Length[A].
However, calculation of Inverse is much more faster than the calculation of PseudoInverse. For the
resistance calculations, it will be enough to compute Inverse[A + 1/k] or Inverse[A - 1/k].
When we do this for Tau calculations, we consider adding (Tr[B] + 1)/k rather
than (Tr[B])/k.*)

Valence[A_] :=
    (Count[A[[#]], Except[0]]-1)&

Genus[A_] :=
    Count[A, Except[0], 2]/2 - 2 Length[A] + 1

CompleteGraphM[k_] :=
    Module[ {kk = (k (k - 1))/2},
        1/2 (-1 + k) k^2 IdentityMatrix[k] - kk
    ]

TotalLength[A_] :=
    Total[1/DeleteCases[DeleteCases[UpperTriangularize[-A, 1], 0, 2],0.,2], 2]

ResistanceMatrix[A_] :=
    Block[ {k, B},
        k = Length[A];
        B = Inverse[A + 1/k];
        Table[
        B[[i, i]] + B[[j, j]] - 2 B[[i, j]], {i, 1, k}, {j, 1, k}]
    ]

AdjacentToLaplacian[A_] :=
    DiagonalMatrix[-Total[A]] + A

ThetaConstant[A_, Q_] :=
    Block[ {k, B, vl},
        k = Length[A];
        B = Inverse[A + 1/k];
        vl = Valence[A][#]&/@ Range[k];
        If[ Q === 0,
            2 Total[Table[(vl[[p]] - 2 ) (vl[[q]] - 2 ) (B[[p, p]] + B[[q, q]] - 2 B[[p, q]]), {p, 1,
               k}, {q, p, k}], 2],
            2 Total[Table[(vl[[p]] - 2 + 2 Q[[p]]) (vl[[q]] - 2 +
            2 Q[[q]]) (B[[p, p]] + B[[q, q]] - 2 B[[p, q]]), {p, 1,k}, {q, p, k}], 2]
        ]
    ]

TauConstant[A_] :=
    Block[ {k, B, S, f},
        k = Length[A];
        B = Inverse[A + 1/k];
        f[a_, {b_, c_}] :=
            (
                If[ A[[b, c]] =!= 0 && A[[b, c]] =!= 0. && b > c,
                    S = (B[[c, c]] + B[[b, b]] - 2 B[[c, b]]);
                    -1/
                              12 (1/A[[b, c]] + S)^2 *A[[b, c]] -
                          A[[b, c]]/4 (B[[c, c]] - B[[b, b]])^2,
                    0
                ]);
        Total[MapIndexed[f, A, {2}], 2] + (Tr[B] - 1)/k
    ]

PhiConstant[A_, Q_] :=
    Block[ {k, B, vl, g, f, theta, tau},
        k = Length[A];
        B = Inverse[A + 1/k];
        vl = Valence[A][#] & /@ Range[k];
        g = If[ Q === 0,
                Total[vl]/2 - k + 1,
                Total[vl]/2 - k + 1 + Total[Q]
            ];
        theta =
         If[ Q === 0,
             2 Total[Table[(vl[[p]] - 2) (vl[[q]] - 2) (B[[p, p]] +
                   B[[q, q]] - 2 B[[p, q]]), {p, 1, k}, {q, p, k}], 2],
             2 Total[Table[(vl[[p]] - 2 + 2 Q[[p]]) (vl[[q]] - 2 +
                   2 Q[[q]]) (B[[p, p]] + B[[q, q]] - 2 B[[p, q]]), {p, 1,
                 k}, {q, p, k}], 2]
         ];
        f[a_, {b_, c_}] :=
            (If[ A[[b, c]] =!= 0 && A[[b, c]] =!= 0. && b > c,
                 S = (B[[c, c]] + B[[b, b]] - 2 B[[c, b]]);
                 -1/12 (1/A[[b, c]] + S)^2*A[[b, c]] -
                  A[[b, c]]/4 (B[[c, c]] - B[[b, b]])^2,
                 0
             ]);
        tau = Total[MapIndexed[f, A, {2}], 2] + (Tr[B] - 1)/k;
        (5 g - 2)/g tau + theta/(4 g ) - TotalLength[A]/4
    ]


End[]

EndPackage[]
\end{verbatim}
}

\end{document}